# A Computationally Efficient Vectorized Implementation of Localizing Gradient Damage Method in MATLAB

Subrato Sarkar[*]

Centre for Modeling and Simulation in Medicine,
Rensselaer Polytechnic Institute, Troy, New York, USA, 12180
*Corresponding author e-mail: ssarkar@me.iitr.ac.in; sarkas6@rpi.edu

## ABSTRACT

In this work, a recently developed fracture modeling method called localizing gradient damage method (LGDM) is implemented in MATLAB. MATLAB is well-known in the computational research community for its simple and easy-to-learn coding interface. As a result, MATLAB is generally preferred for the initial development (prototyping) of computational models by researchers. However, MATLAB-developed codes are seldom used for large-scale simulations (after initial development is complete) due to their computational inefficiency. Hence, a computationally efficient implementation of LGDM using MATLAB vectorization is presented in this work. The choice of LGDM (as the fracture modeling method) is based on its thermodynamically consistent formulation built upon the micromorphic framework. Moreover, the non-linear coupled field formulation of LGDM makes it suitable for testing the computational efficiency of vectorized MATLAB implementation in a non-linear finite element setting. It is shown in this work that the vectorized MATLAB implementation can save significant computational resources and time as compared to non-vectorized implementations (that are parallelized with MATLAB `parfor`). The vectorized MATLAB implementation is tested by solving numerical problems in 1D, 2D and 3D on a consumer-grade PC, demonstrating the capability of vectorized implementation to run simulations efficiently on systems with limited resources. The sample source codes are provided as supplementary materials that would be helpful to researchers working on similar coupled field models.

**Keywords**: MATLAB; FEM; localizing; Gradient damage; Coupled Field; Nonlinear

## 1. INTRODUCTION

Computational modeling of fracture in materials is a fairly established field of research nowadays, with applications ranging from nano/micro scale to large-scale structures (Brunner, 2020; Mohammadnejad *et al*., 2018; Doblare *et al*., 2004). However, in computational fracture modeling, the learning curve is often very steep for new researchers, and even experienced researchers (looking for an application of fracture modeling in other fields) find it hard to get started. This is because a substantial effort is needed to develop these computational fracture models from scratch and also the non-availability of open-source codes. Fortunately, nowadays, many researchers are making efforts to publish/share codes of computational fracture modeling methods. This would undoubtedly lead to faster development and application of computational fracture models to new frontier areas. The positive effect created by the availability of open-source codes is evident from recent rapid developments in artificial intelligence research (Engler, 2022).

In computational fracture modeling, numerous works have been published recently on computer implementations using different software. The most popular software used in these works are Abaqus, COMSOL and MATLAB (Abaqus, 2014; COMSOL, 2022; MATLAB, 2022). Among these software, Abaqus and COMSOL are commercial finite element (FE) packages, while MATLAB is a general-purpose programming package. It is emphasized that Abaqus and COMSOL packages are based on the finite element method (FEM) and are only capable of modeling problems formulated using FEM (or some variant of it e.g., extended FEM). However, the implementations in MATLAB are not limited to FEM and are more versatile. Various works on fracture modeling methods implemented in Abaqus, COMSOL and MATLAB are discussed in the following.

Computational fracture modeling is broadly categorized into discrete and smeared models (Ambati *et al*., 2014). In the discrete models, the crack is defined explicitly using either conformal mesh or nodal enrichment (Borst *et al*., 2004; Jha *et al*., 2022). In comparison, the smeared methods have diffused definitions of crack (i.e. phase field models (Patil *et al*., 2018) and gradient damage models (Sarkar *et al*., 2020a, b; Sarkar *et al*., 2021a; Sarkar *et al*., 2022a; Bansal *et al*., 2022). Notable works on the computer implementation of discrete methods include extended FEM (XFEM) (Giner *et al*., 2009; Shi *et al*., 2010; Rokhi and Shariati, 2013; Sutula *et al*., 2018; Ding *et al*., 2020; Jafari *et al*., 2022) and cohesive zone method (Park and Paulino, 2012). For smeared methods, the popular works on computer implementation include the phase field method (Msekh *et al*., 2015; Molnar and Gravouil, 2017; Zhou *et al*., 2018; Fang *et al*., 2019; Molnar *et al*., 2020; Chen and Wu, 2022; Rahaman, 2022), peridynamics (Huang et al., 2019; Bie *et al*., 2020, Jenabidehkordi and Rabczuk, 2022) and gradient damage method (Seupel *et al*., 2018; Sarkar *et al*., 2019a, Sarkar *et al*., 2019b, Sarkar *et al*., 2022b; Zhang *et al*., 2022).

Despite the abundance of implementations in Abaqus and COMSOL, it is generally found that the computer implementations in MATLAB are easy to understand and modify. The fundamental reasons for this are elucidated in Figure 1. However, it is acknowledged that simulations in MATLAB are usually slower and more inefficient than commercial FE codes like Abaqus and COMSOL. Perhaps, due to MATLAB's computational inefficiency, it is mainly used for prototyping (early development and understanding) a new model/idea that is later implemented in Abaqus or COMSOL for computational efficiency. Hence, it is shown in this work that a MATLAB implementation can be significantly more efficient for simulating non-linear phenomena like a fracture. Besides computational efficiency, understandability and modifiability are added benefits of the presented MATLAB implementations.

| **MATLAB Implementations** | **Abaqus/COMSOL Implementations** |
|---|---|
| 1. *Easily Understandable*:<br>　○ Normally entire FE code available<br>　○ All aspects of implementation understood from code (Elemental computations, assembly, variable update etc.)<br>　○ Simple syntax (interpreted language), no coding specific knowledge needed<br><br>2. *Easy to Setup and Modify*:<br>　○ Simply open MATLAB and run the code<br>　○ Modifications carried out on the fly (interpreted language)<br><br>3. *Generally Slow and Inefficient* | 1. *Not Easily Understandable*:<br>　○ Code available as subroutine (Abaqus) or no code available (COMSOL)<br>　○ All aspects of implementation not clear form the available code or files (e.g. Elemental computations, variable update)<br>　○ Fortran (in Abaqus) needs coding specific knowledge and COMSOL needs interface specific knowledge<br><br>2. *Not Easy to Setup and Modify*:<br>　○ Abaqus: Linking of compiler and solver needed to run and modify subroutines<br>　○ COMSOL: Full understanding of the user interface needed to run and modify<br><br>3. *Generally Fast and Memory Efficient* |

**Figure 1**: Major differences in the available open-source MATLAB and Abaqus/COMSOL Implementations

In the present work, a recently developed non-linear fracture modeling method called localizing gradient damage method (LGDM) is used to model fracture/damage. The LGDM is well-established as an accurate and efficient method of fracture modeling (Poh and Sun, 2016; Huang *et al.*, 2022). In this regard, the novelties of the present work are as follows,

- A MATLAB implementation of LGDM is presented in 1D, 2D and 3D. The presented implementation includes both non-vectorized and vectorized MATLAB codes.
- A comparison between non-vectorized and vectorized implementations shows significant savings in computation time and memory usage by the vectorized implementation.
- The differences in the structure of a non-vectorized and vectorized MATLAB code are discussed in detail, which can be used to vectorize any coupled-field non-linear method.

- The open-source codes for the implementation are shared for the benefit of the research community.

This paper is structured into four sections. Section 1 discusses the introduction, motivation and novelties of the current study. Then, the LGDM is briefly reviewed in Section 2, along with a discussion on the MATLAB implementation aspects. A comparative study of non-vectorized and vectorized MATLAB implementation is presented in the Results and Discussion (Section 3). Finally, the major outcomes of the present study are highlighted in the Conclusion (Section 4).

## 2. LGDM AND ITS IMPLEMENTATION ASPECTS

This section briefly reviews the formulation of localizing gradient damage method (LGDM) in the first sub-section. The details of non-vectorized and vectorized MATLAB implementation of LGDM are described in the following sub-sections.

### 2.1 Localizing Gradient Damage Method (LGDM)

The localizing gradient damage method (LGDM) is a thermodynamically consistent fracture modeling method and has also been proven to be more accurate and computationally efficient (Huang *et al*., 2022). The accuracy of LGDM can be attributed to the micromorphic framework, in which, a morphic variable is introduced to account for the fracture processes at the underlying micro-continuum. This morphic variable is introduced in addition to the traditional kinematic variables (Poh and Sun, 2016). The micromorphic framework enables LGDM to include the effects of fluctuating micro-level responses during a fracture that are otherwise neglected in the macroscopic continuum theory.

In the micromorphic framework, the free energy density is assumed such that it includes energy due to higher-order stresses ($\bar{\sigma}$ and $\bar{\bar{\xi}}$) associated with the micromorphic variable ($\bar{\varepsilon}_{eq}$) and its gradient ($\nabla \bar{\varepsilon}_{eq}$). This additional energy due to higher order terms is in addition to the standard strain energy. Hence, the free energy density ($\Psi$) is expressed as,

$$\Psi = \frac{1}{2}(1-D)\boldsymbol{\varepsilon} : {}^4\mathbf{C} : \boldsymbol{\varepsilon} + \frac{1}{2}h\left(\varepsilon_{eq} - \bar{\varepsilon}_{eq}\right)^2 + \frac{1}{2}ghc\left(\nabla \bar{\varepsilon}_{eq} \cdot \nabla \bar{\varepsilon}_{eq}\right) \qquad (1)$$

In Eq. (1), the first term denotes the standard elastic strain energy characterized by the fourth-order elasticity tensor (${}^4\mathbf{C}$), local/macro tensorial strain ($\boldsymbol{\varepsilon}$) and a scalar damage variable (*D*). The subsequent second and third terms are non-standard that characterize the coupling interactions (macro-micro interactions) and micro-micro interactions (at micro-scale) occurring in the fracture process zone. Note that the micromorphic variable is called micro-equivalent strain ($\bar{\varepsilon}_{eq}$) in this paper. The difference between macro-equivalent and micro-equivalent strains in the second term ($\varepsilon_{eq} - \bar{\varepsilon}_{eq}$) accounts for coupling (micro-macro) interactions whose magnitude is quantified by a parameter called the coupling modulus (*h*). The contribution of micro-micro interactions occurring at the micro-continuum is included through the gradient of the micro-equivalent strains ($\nabla \bar{\varepsilon}_{eq}$) in the third term of Eq. (1). Apart from the coupling modulus (*h*), a couple of additional parameters are defined in the third term, i.e. the gradient parameter (*c*) and the interaction parameter (*g*). The significance of these parameters is discussed extensively in the author's previous work (Sarkar *et al.*, 2019a, Sarkar *et al.*, 2020b) and thus avoided here for brevity. However, the expressions for the damage law (*D*), macro equivalent strain ($\varepsilon_{eq}$) and interaction function (*g*) used in the present work are defined in Appendix A.

From Eq. (1), the constitutive relations are obtained by following the Coleman-Noll procedure on the free energy density function ($\Psi$) as (Poh and Sun, 2016),

$$\boldsymbol{\sigma} = (1-D)\,^4\mathbf{C} : \boldsymbol{\varepsilon} + h(\varepsilon_{eq} - \bar{\varepsilon}_{eq})\frac{\partial \varepsilon_{eq}}{\partial \boldsymbol{\varepsilon}} \qquad (2)$$

$$\bar{\sigma} = h(\varepsilon_{eq} - \bar{\varepsilon}_{eq}) \qquad (3)$$

$$\bar{\boldsymbol{\xi}} = ghc(\nabla \bar{\varepsilon}_{eq}) \qquad (4)$$

Further, the substitution of constitutive equations into the energy balance yields the following governing equations (Poh and Sun, 2016),

$$\nabla \cdot \boldsymbol{\sigma} = \mathbf{0} \text{ in domain } \Omega \qquad (5)$$

$$\bar{\sigma} = \nabla \cdot \bar{\boldsymbol{\xi}} \text{ in domain } \Omega \qquad (6)$$

with boundary conditions,

$$\boldsymbol{\sigma} \cdot \mathbf{n} = \mathbf{t} \text{ on boundary } \partial\Omega \qquad (7)$$

$$\bar{\boldsymbol{\xi}} \cdot \mathbf{n} = \zeta \text{ on boundary } \partial\Omega \qquad (8)$$

where **t** and $\zeta$ are surface traction and higher-order traction, respectively. Using the method of weighted residuals and suitable finite element discretization, the governing equations in Eqs. (5)-(8) can be written after consistent linearization (Sarkar *et al.*, 2020b),

$$\int_\Omega \mathbf{B}_u^T \delta\boldsymbol{\sigma} \, d\Omega = \int_{\partial\Omega} \mathbf{N}_u^T \mathbf{t} \, d\partial\Omega - \int_\Omega \mathbf{B}_u^T \boldsymbol{\sigma}_{i-1} \, d\Omega \qquad (9)$$

$$\int_\Omega \mathbf{N}_{\bar{\varepsilon}}^T \delta\bar{\sigma} \, d\Omega + \int_\Omega \mathbf{B}_{\bar{\varepsilon}}^T \delta\bar{\boldsymbol{\xi}} \, d\Omega = \int_\Omega \mathbf{N}_{\bar{\varepsilon}}^T \bar{\sigma}_{i-1} \, d\Omega + \int_\Omega \mathbf{B}_{\bar{\varepsilon}}^T \bar{\boldsymbol{\xi}}_{i-1} \, d\Omega \qquad (10)$$

where the symbol δ denotes a linearized increment while **N** and **B** are the shape functions and derivatives associated with a variable denoted by an appropriate subscript. For discretization, the shape functions for displacement ($u$) are taken as quadratic, and for micro-equivalent strain ($\bar{\varepsilon}_{eq}$) are taken as linear. Eqs. (9)-(10) can be expressed in compact matrix form after the substitution of linearized variables as,

$$\begin{bmatrix} \mathbf{K}^{uu}_{i-1} & \mathbf{K}^{u\bar{\varepsilon}}_{i-1} \\ \mathbf{K}^{\bar{\varepsilon}u}_{i-1} & \mathbf{K}^{\bar{\varepsilon}\bar{\varepsilon}}_{i-1} \end{bmatrix} \begin{Bmatrix} \delta\widetilde{\mathbf{u}} \\ \delta\widetilde{\bar{\boldsymbol{\varepsilon}}}_{eq} \end{Bmatrix} = \begin{bmatrix} \mathbf{F}^{u}_{i-1} \\ \mathbf{F}^{\bar{\varepsilon}}_{i-1} \end{bmatrix} \tag{11}$$

where,

$$\mathbf{K}^{uu}_{i-1} = \int_\Omega \mathbf{B}_u^T (1 - D_{i-1})\,^4\mathbf{C}\,\mathbf{B}_u\, d\Omega \tag{11a}$$

$$\mathbf{K}^{u\bar{\varepsilon}}_{i-1} = -\int_\Omega \mathbf{B}_u^T \left\{ {}^4\mathbf{C}\boldsymbol{\varepsilon}_{i-1} \left[\frac{\partial D}{\partial \bar{\kappa}}\right]_{i-1} \left[\frac{\partial \bar{\kappa}}{\partial \bar{\varepsilon}_{eq}}\right]_{i-1} + h \left[\frac{\partial \varepsilon_{eq}}{\partial \boldsymbol{\varepsilon}}\right]_{i-1} \right\} \mathbf{N}_{\bar{\varepsilon}}\, d\Omega \tag{11b}$$

$$\mathbf{K}^{\bar{\varepsilon}u}_{i-1} = -\int_\Omega \mathbf{N}_{\bar{\varepsilon}}^T h \left[\frac{\partial \varepsilon_{eq}}{\partial \boldsymbol{\varepsilon}}\right]_{i-1} \mathbf{B}_u\, d\Omega \tag{11c}$$

$$\mathbf{K}^{\bar{\varepsilon}\bar{\varepsilon}}_{i-1} = \int_\Omega \left\{ \left( \mathbf{N}_{\bar{\varepsilon}}^T h + \mathbf{B}_{\bar{\varepsilon}}^T h c \left[\frac{\partial g}{\partial \bar{\kappa}}\right]_{i-1} \nabla \bar{\varepsilon}_{eq_{i-1}} \right) \mathbf{N}_{\bar{\varepsilon}} + \mathbf{B}_{\bar{\varepsilon}}^T g h c \mathbf{B}_{\bar{\varepsilon}} \right\} d\Omega \tag{11d}$$

$$\mathbf{F}^{u}_{i-1} = \int_{\partial\Omega} \mathbf{N}_u^T \mathbf{t}\, d\partial\Omega - \int_\Omega \mathbf{B}_u^T \boldsymbol{\sigma}_{i-1}\, d\Omega \tag{11e}$$

$$\mathbf{F}^{\bar{\varepsilon}}_{i-1} = \int_\Omega \mathbf{N}_{\bar{\varepsilon}}^T h (\varepsilon_{eq} - \bar{\varepsilon}_{eq})_{i-1}\, d\Omega + \int_\Omega \mathbf{B}_{\bar{\varepsilon}}^T g h c (\nabla \bar{\varepsilon}_{eq})_{i-1}\, d\Omega \tag{11f}$$

The parameters appearing in Eqs. 11(a-f) are outlined in Appendix A for reference. A non-linear solution procedure using Newton's method is adopted for solving the system of equations in Eq. (11). The incremental-iterative solution procedure is used. The following sub-sections outline the MATLAB implementation methodology adopted for simulating fracture using the aforementioned formulation.

**2.2 MATLAB Implementation of LGDM**

This sub-section describes the implementation of LGDM in MATLAB. The non-vectorized implementation is discussed first, followed by the vectorized implementation. The discussion is elaborated using screenshots of MATLAB code, making it easier to understand. 2D implementation is used to describe the MATLAB code. The 2D implementation is chosen for discussion because it is simpler to understand than the 3D implementation and is not oversimplified like the 1D implementation. For brevity, the discussion is carried out only on those parts of the algorithm (and code) that are significantly different in the non-vectorized and vectorized MATLAB code. However, the interested reader is referred to the MATLAB codes shared as supplementary materials for a detailed look.

## 2.2.1 Non-vectorized Implementation

In this sub-section, the non-vectorized MATLAB implementation of LGDM is discussed. A code implementation that majorly uses for-loops is called a non-vectorized implementation (MATLAB, 2022). A non-vectorized implementation is generally adopted during the initial stages of development (while developing in-house codes) because they are easy to understand and modify. Algorithm 1 shows various steps used in the non-vectorized implementation of LGDM. Note that apart from using for-loops for loadsteps (line 2) and iterations (line 3), the for-loops are also used for elemental computations (lines 4-10), assembly (lines 9-11) and variable updates (lines 14-24).

---

**Algorithm 1:** Solution procedure for LGDM using non-vectorized code

**Input:** Geometry, Material properties, Loads and boundary conditions
**Output:** Damage and Structural response

```
1   Initialize: ũ = 0; ε̃̄_eq = 0; SDVs = 0                     // As Cell Arrays
2   for n ← 1 to N do                                          // Total loadsteps = N
3     for i ← 1 to I do                                        // Total Iterations = I
4       for el ← 1 to nel do                                   // Total Elements = nel
5         for igp ← 1 to ngp do                                // Elemental GPs = ngp
6           Compute K_local: K^uu_{i-1}; K^uε̄_{i-1}; K^ε̄u_{i-1}; K^ε̄ε̄_{i-1}   ▶ Eqs. 11a-d
7           Compute F_local: F^u_{i-1}; F^ε̄_{i-1}              ▶ Eqs. 11e-f
8         end
9         Unroll: K_local & F_local
10      end
11      Assemble: K and F                                      // Assembly
12      Solve: Kδũ = F                                         ▶ Eq. 11
13      ũ_i = ũ_{i-1} + δũ; ε̃̄_eq_i = ε̃̄_eq_{i-1} + δε̃̄_eq     // Update Primary Var
14      for el ← 1 to nel do                                   // Update SDVs
15        for igp ← 1 to ngp do
16          Update:
17            ε           (Strain Tensor)
18            ε_eq        (Local Equivalent Strain)            ▶ Eq. A.3
19            ∂ε_eq/∂ε    (Derivatives of Equivalent Strain)   ▶ Eq. A.3
20            g           (Interaction Function)               ▶ Eq. A.4
21            D           (Damage)                             ▶ Eq. A.1
22            σ           (Stress Tensor)                      ▶ Eq. 2
23        end
24      end
25    end
26    Check Convergence: ‖δũ/ũ‖ < Tol; ‖δε̃̄_eq/ε̃̄_eq‖ < Tol
27    If Converged GOTO Line 2
28    Else GOTO Line 3
29  End
```

---

Specifically, the present discussion focuses on two typical operations during each iteration of a non-linear finite element simulation, which are,

(a) Elemental (local) computations of **K** & **F** and their assembly (Algorithm 1, lines 4-11)

(b) Variable updates after obtaining a solution (Algorithm 1, lines 13-24)

These operations in non-vectorized implementations are usually carried out using for-loops. These for-loops, if inefficiently handled, can incur most of the computational cost. The discussion on the abovementioned two operations is as follows,

(a) *Elemental Computations and Assembly*

In Algorithm 1, the elemental (local) computations of stiffness matrix ($\mathbf{K}_{local}$) and force vector ($\mathbf{F}_{local}$) are carried out between lines 4-10, and their assembly is in lines 9 & 11. These elemental computations and assembly are shown through a MATLAB code snippet in Figure 2. Note the use of for-loop over elements to compute elemental (local) contributions to the global stiffness ($\mathbf{K}$) and force vector ($\mathbf{F}$). The code in Figure 2 has three distinct sections involving computations of shape function derivatives/Jacobian, computation of elemental $\mathbf{K}$ & $\mathbf{F}$ and assembly.

It is pointed out that a couple of unique features available in MATLAB are used to optimize this non-vectorized implementation, which are **(i)** Parallel for-loop (`parfor`) and **(ii)** `cell` arrays. The `parfor` is a MATLAB feature that parallelizes the for-loop (distributes for-loop computations) across multiple CPU processor cores. It is later shown in the numerical problems that using `parfor` (parallel for-loop) reduces computation time in 2D and 3D simulations. Besides, it can be observed in MATLAB codes (Figures 2, 3, 4 & 5) that a different type of array called `cell` array is used. The `cell` arrays are designed to store arrays within them (i.e. arrays within an array) and are indexed using curly braces '{}'. In addition to increasing the readability of code (by avoiding multidimensional arrays), the `cell` arrays are preferred by MATLAB for executing `parfor` loops.

The different sections of the code shown in Figure 2 are discussed below.

- *Computations of shape function derivatives/Jacobian*: The computations of `B_mat_u` and `J` representing shape function derivatives and element Jacobian (Figure 2, lines 3-4) used for Gauss integration are elaborated through a MATLAB code snippet in Figure 3. Note the use of for-loop for each Gauss point of an element. In Figure 3, the sizes of all the arrays used in the computation are mentioned in subscripted square brackets.

- *Computations of elemental $\mathbf{K}$ and $\mathbf{F}$*: The different components of $\mathbf{K}_{local}$ and $\mathbf{F}_{local}$ (Algorithm 1, lines 6-7) are computed using function subroutines (Figure 2, lines 7-18). Within these function subroutines, the contributions of individual Gauss points are computed through for-loops. These computations at individual Gauss points are shown for $K_{i-1}^{uu}$ (1$^{st}$ component of **K**$_{local}$) in Figure 4 as a function subroutine. Similar function subroutines with Gauss point for-loops are used for other components of **K**$_{local}$ and **F**$_{local}$. The Gauss point for-loops within these function subroutines run sequentially and cannot be efficiently parallelized. In other words, the `parfor` can only parallelize the overall elemental computations (Figure 2, line 2), and computations within each element, i.e. at the Gauss points (Figure 2, lines 7-18), must be sequential. This is a limitation of the non-vectorized implementation that leads to computation overhead.

- *Assembly of K and F*: The elemental computations for all the components of **K**$_{local}$ and **F**$_{local}$ are saved in local arrays called `k_local` and `f_local` (in Figure 2, lines 20-21). These local arrays are then unrolled into vectors (Figure 2, lines 24-25) and used to create `sparse` matrices (Figure 2, lines 29-30). Note the extensive use of `cell` arrays to unroll and assemble elemental matrices.

```
1   % Build elemental contributions to K and F..........Start element Loop
2   parfor el = 1:nel
3       [B_mat, B_mat_tilde, J, ~, N_mat_tilde] = gp_data(gp, ...
4           c_mat_u(el,:), B, ngp, nen_u, nen_e);
5
6       % Stiffness Matrix
7       [k_uu] = Build_kuu(B_mat, D{el,1}, D_mat, J, gw, ngp, nen_u);
8       [k_ue, dD_dk] = Build_kue(B_mat, D_mat, strain{el,1},...
9           N_mat_tilde, k0{el,1}, k_gp{el,1}, ki, J, gw, ngp, nen_u, ...
10          nen_e, beta, alpha);
11      [k_eu] = Build_keu(B_mat, N_mat_tilde, strain{el,1}, ...
12          J, gw, ngp, nen_u, nen_e, k, c1, c2, c3, h);
13      [k_ee,G] = Build_kee(B_mat_tilde, N_mat_tilde, c, ...
14          J, gw, ngp, nen_e, e_cell{el,1}, D{el,1}, h, R, nta, dD_dk);
15
16      % RHS
17      [f_e] = Build_fee(B_mat_tilde, N_mat_tilde, loc_e_eqv{el,1},...
18          gp_e_eqv{el,1}, e_cell{el,1}, G, J, gw, ngp, nen_e, h, c);
19      f_u = -k_uu*u_cell{el,1};
20      k_local = [k_uu k_ue; k_eu k_ee];
21      f_local = [f_u; f_e];
22
23      % Building row, col, value for assembly
24      [krow{el,1},kcol{el,1},kval{el,1},frow{el,1},fval{el,1}] = ...
25          Build_sparse(connect_u{el,1},connect_e{el,1}, k_local', f_local);
26  end
27
28  % Assemble K, F
29  K = sparse(cat(1,krow{:}),cat(1,kcol{:}),cat(1,kval{:}));
30  F = sparse(cat(1,frow{:}),1,cat(1,fval{:}));
```

Annotations:
- For-Loop over each Element
- Shape Fn Derivatives and Jacobian (Figure 3)
- **K** and **F** Computation Subroutines (Figure 4)
- **K** and **F** Assembly

**Figure 2**: Code snippet showing the elemental computations and assembly

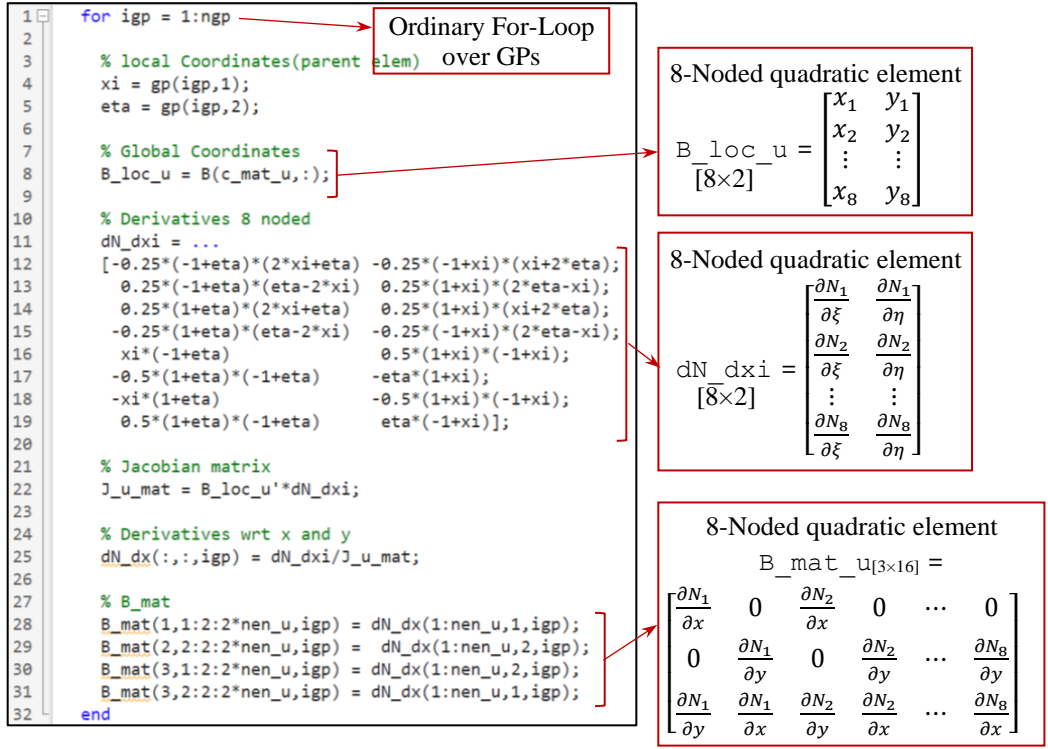

**Figure 3**: Code snippet showing the computation of strain displacement matrix (`B_mat_u`) and Jacobian (`J`)

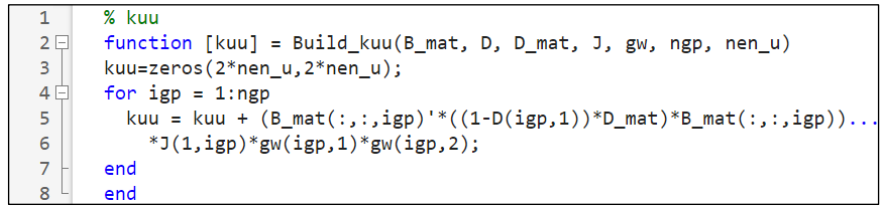

**Figure 4**: Code snippet showing the elemental computation function subroutine for $\mathbf{K}_{i-1}^{uu}$

(b) *Variables Update*

The primary and solution-dependent variables (SDVs) are updated (Algorithm 1, lines 13-24) after the incremental values of the solution ($\delta \widetilde{\mathbf{u}}$ and $\delta \widetilde{\bar{\boldsymbol{\varepsilon}}}_{eq}$) are obtained by solving the system of equations. The primary variables are updated at the nodal points, while SDVs are updated at the Gauss points. A MATLAB code snippet in Figure 5a shows the primary variable and SDV updates. The primary variables are updated using the incremental solution vector (Figure 5a, line 2), and the SDVs are updated using an elemental for-loop (Figure 5a, lines 9-16). It is noted that the `parfor` is again used to parallelize the elemental for-loop in the SDV update (Figure 5a, line 9).

The function subroutine used for SDV update is shown in Figure 5b. In this function subroutine, the Gauss point variables are updated using a for-loop over Gauss points, similar to the stiffness matrices mentioned previously. Consequently, the for-loop within the function subroutine runs sequentially, leading to increased computational costs.

```matlab
1   % Update primary variables
2   u = u + du;
3   for el = 1:nel
4       u_cell{el,1} = u(connect_u{el,1},1);
5       e_cell{el,1} = u(connect_e{el,1},1);
6   end
7
8   % Update SDVs
9   parfor el = 1:nel
10      [B_mat, ~, ~, ~, N_mat_tilde] = gp_data(gp, ...
11          c_mat_u(el,:), B, ngp, nen_u, nen_e);
12      [strain{el,1}, loc_e_eqv{el,1}, gp_e_eqv{el,1}, k_gp{el,1}, ...
13          D{el,1}] = ...
14          update_variables(k0{el,1}, ki, u_cell{el,1}, e_cell{el,1}, ...
15          B_mat, N_mat_tilde, ngp, k, alpha, beta, c1, c2, c3);
16  end
```
(Lines 3–6 annotated: For-Loop over Elements)

**(a)** Update variables

```matlab
1   % Update variables.....................................................
2   function [strain, loc_e_eqv, gp_e_eqv, k_gp, D] = ...
3       update_variables(k0, ki, u, e, B_mat, N_mat_tilde, ngp, k, ...
4       alpha, beta, c1, c2, c3)
5
6   strain = zeros(ngp,3);    % Strains
7   loc_e_eqv = zeros(ngp,1); % local equivalent strain
8   gp_e_eqv = zeros(ngp,1);  % equivalent strains at GP from solution
9   k_gp = zeros(ngp,1);      % History strain at GP
10  D = zeros(ngp,1);         % Damage
11  for igp = 1:ngp
12      strain(igp,:) = B_mat(:,:,igp)*u;
13      loc_e_eqv(igp,:) = c1*(strain(igp,1) + strain(igp,2))+(1/(2*k))*...
14          ((c2*(strain(igp,1) + strain(igp,2)))^2 + ...
15          (c3*(strain(igp,1)^2 + strain(igp,2)^2 - ...
16          strain(igp,1)*strain(igp,2) + 3*(strain(igp,3)^2))))^0.5;
17      gp_e_eqv(igp,:) = N_mat_tilde(:,igp)'*e;
18
19      if (gp_e_eqv(igp,:) < k0(igp,:))
20          k_gp(igp,:) = k0(igp,:);
21      else
22          k_gp(igp,:) = gp_e_eqv(igp,:);
23      end
24      if (k_gp(igp,:) < ki)
25          D(igp,:) = 0.0;
26      else
27          D(igp,:) = 1-(ki/k_gp(igp,:))*(1-alpha+alpha*...
28              exp(-beta*(k_gp(igp,:) - ki)));
29      end
30  end
31  end
```
(Line 11 annotated: For-Loop over GPs; Lines 19–29 annotated: If-Else conditions for consistent variable updates)

**(b)** Update variables function subroutine

**Figure 5**: Code snippets for update variables and its function subroutine

### 2.2.2 Vectorized Implementation

The non-vectorized parts (that use for-loops) of the MATLAB code mentioned in the previous sub-section are converted to a vectorized code and described in this sub-section. The conversion from non-vectorized to vectorized code means eliminating the for-loops with vectorized operations. In other words, the operations previously carried out using for-loops are converted to entry-wise operations on arrays that do not need for-loops. The entry-wise operations are faster in MATLAB because the entire array is operated upon simultaneously instead of each entry through a for-loop. These operations are carried out using special commands in MATLAB. The commonly used special MATLAB commands in current vectorization are `repmat`, `kron`, `reshape`, `sum`, entry-wise multiplication (`.*`) and entry-wise division (`./`) (MATLAB, 2022). An 'array entry' is also usually called an 'element' of the array; however, the term 'element' is reserved in this work for finite elements. For example, an entry at $i^{th}$ row and $j^{th}$ column of array **A** is denoted as $\mathbf{A}(i,j)$. The solution procedure for LGDM after MATLAB vectorization is shown in Algorithm 2.

---

**Algorithm 2:** Solution procedure for LGDM using vectorized code

**Input:** Geometry, Material properties, Loads and boundary conditions
**Output:** Damage and Structural response

| | | |
|---|---|---|
| 1 | **Initialize:** $\tilde{\mathbf{u}} = \mathbf{0}$; $\tilde{\bar{\boldsymbol{\varepsilon}}}_{eq} = \mathbf{0}$; SDVs = **0** | // As Vectors |
| 2 | **for** $n \leftarrow 1$ **to** $N$ **do** | // Total loadsteps = $N$ |
| 3 |    **for** $i \leftarrow 1$ **to** $I$ **do** | // Total Iterations = $I$ |
| 6 |       **Compute and Assemble K** | ▶ Eqs. 11a-d |
| 7 |       **Compute and Assemble F** | ▶ Eqs. 11e-f |
| 8 |       **Solve:** $\mathbf{K}\delta\tilde{\mathbf{u}} = \mathbf{F}$ | ▶ Eq. 11 |
| 9 |       $\tilde{\mathbf{u}}_i = \tilde{\mathbf{u}}_{i-1} + \delta\tilde{\mathbf{u}}$; $\tilde{\bar{\boldsymbol{\varepsilon}}}_{eq_i} = \tilde{\bar{\boldsymbol{\varepsilon}}}_{eq_{i-1}} + \delta\tilde{\bar{\boldsymbol{\varepsilon}}}_{eq}$ | // Update Primary Var |
| 10 |       **Update:** | // Update SDVs |
| 11 |         $\boldsymbol{\varepsilon}$ (Strain Tensor) | |
| 12 |         $\varepsilon_{eq}$ (Local Equivalent Strain) | ▶ Eq. A.3 |
| 13 |         $\frac{\partial \varepsilon_{eq}}{\partial \boldsymbol{\varepsilon}}$ (Derivatives of Equivalent Strain) | ▶ Eq. A.3 |
| 14 |         $g$ (Interaction Function) | ▶ Eq. A.4 |
| 15 |         $D$ (Damage) | ▶ Eq. A.1 |
| 16 |         $\boldsymbol{\sigma}$ (Stress Tensor) | ▶ Eq. 2 |
| 17 |    **end** | |
| 18 |    **Check Convergence:** $\left\|\frac{\delta\tilde{\mathbf{u}}}{\tilde{\mathbf{u}}}\right\| < \text{Tol}$; $\left\|\frac{\delta\tilde{\bar{\boldsymbol{\varepsilon}}}_{eq}}{\tilde{\bar{\boldsymbol{\varepsilon}}}_{eq}}\right\| < \text{Tol}$ | |
| 19 |    If Converged **GOTO** Line **2** | |
| 20 |    Else **GOTO** Line **3** | |
| 21 | **End** | |

---

After vectorization, the following changes can be observed in Algorithm 2 (vectorized) by comparing it with Algorithm 1 (Non-vectorized),

- For-loops over elements and Gauss points are eliminated in elemental (local) computations (Algorithm 2, Lines 6–7), assembly (Algorithm 2, Lines 6–7) and variables update (Algorithm 2, Lines 10–16).
- The assembly is carried out simultaneously with elemental computations of **K** and **F**, i.e. a separate assembly operation is eliminated.

- MATLAB `parfor` is eliminated because the elemental for-loops are no longer needed. This avoids the additional memory needed to parallelize for-loops through `parfor`.

The initialization of variables (primary and SDVs) in the vectorized implementation is carried out using 1D vectors rather than `cell` arrays. The initializations as 1D vectors lead to easy and efficient handling of the variables by the vectorized code. For example, the damage variable (*D*), which was previously defined as a `cell` array with size `D{nel,1}` in the non-vectorized implementation, is now a vector of size `D(nel*ngp,1)`. Where '{}' is used for indexing a cell array, '()' is used in indexing a standard array/vector, `nel` = number of elements and `ngp` = number of Gauss points in an element.

Similar to the previous sub-section, this discussion is focussed on two typical operations in a non-linear finite element code i.e., the elemental computations with assembly and variables update. The discussion is as follows,

(a) *Elemental Computations and Assembly*

Unlike the non-vectorized implementation, the vectorized implementation is structured such that the elemental computations and assembly are carried out simultaneously. Hence, the compute and assemble operations in lines 6 and 7 of Algorithm 2 are carried out through single-function subroutines for **K** and **F**. These single-function subroutines compute the expressions for the entire model (i.e. for all elements/nodes/Gauss points) at once without explicit for-loops. A couple of important computations in these function subroutines are discussed in the following,

- *Computations of shape function derivatives/Jacobian*: The computation of shape function derivatives (`dN_dx`) and element Jacobian (`J`) used in the Gauss integration are shown through code snippets in Figure 6. The arrays appearing in the computations are shown alongside code expressions. There are mainly two points to note in Figure 6:
  
  o The array structures are such that the expressions are computed for the entire model simultaneously. For example, in lines 31-34, the size of `dN_dx` is [8×(nel*ngp)], which implies that `dN_dx` has values for shape function derivatives corresponding to each node (8) of all the elements (`nel`) for all the GPs (`ngp`).
  
  o The vectorized operations such as entry-wise multiplication (.*)/division (./), `repmat`, `reshape` and `kron` are used on the arrays for loop-free execution.

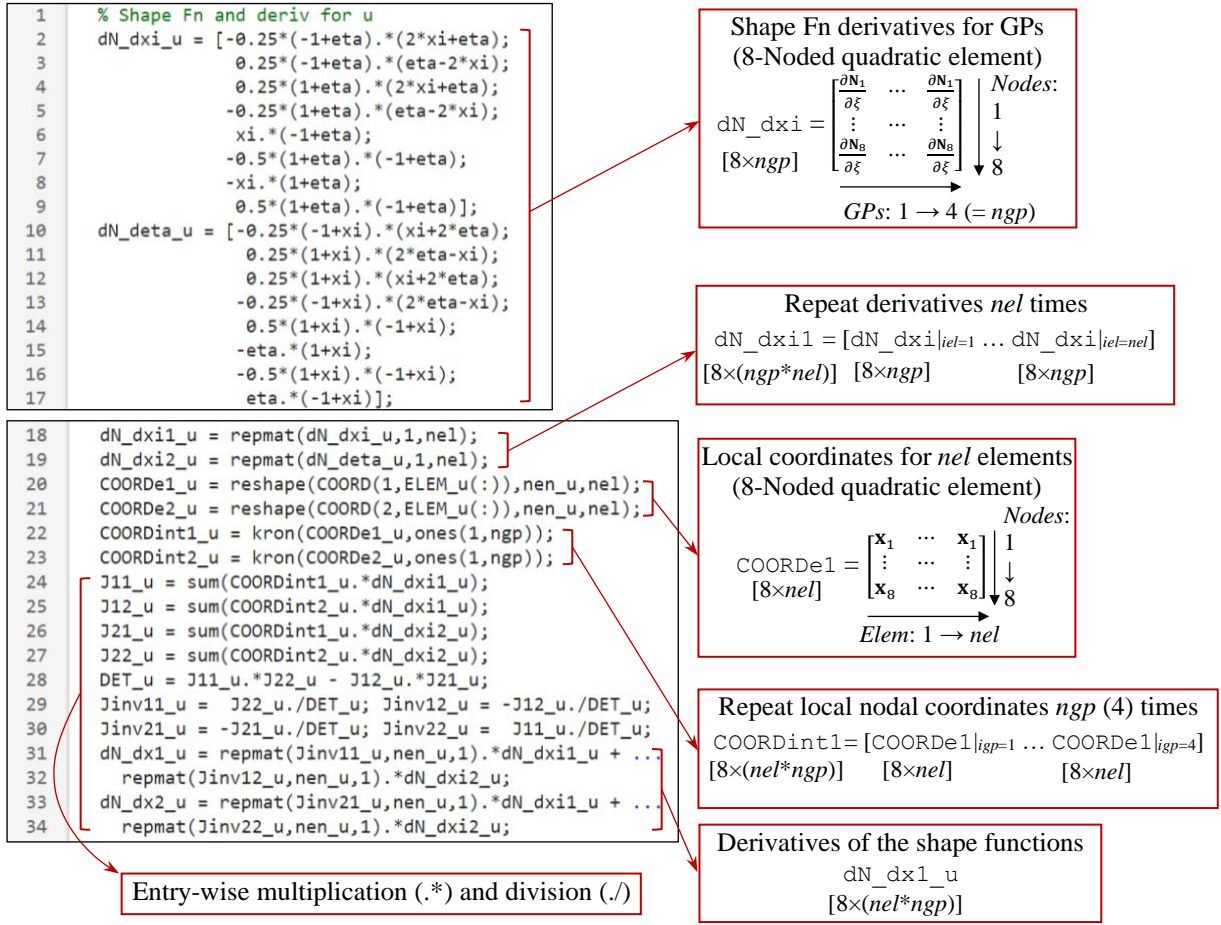

**Figure 6**: Code snippets showing vectorized computation of derivatives of the shape functions (`dN_dx`) and Jacobian (`J`)

- *Computation of global stiffness* (**K**) *matrix*: The stiffness matrix computation for the `Kuu` component is shown through a code snippet in Figure 7. Note that the shape function derivative matrix (`B_mat_u`) and the elasticity matrix (`Duu_mat`) are structured such that the values for all the degrees of freedom (*dofs*) and all Gauss points (*ngp*nel*) are placed at appropriate locations that yield the size of `Kuu` as [*ndof×ndof*]. The computation of `Kuu` at line 25 indicates simultaneous evaluation and assembly. It is emphasized that MATLAB sparse arrays are used extensively to save the memory required and efficiently assemble large arrays (`B_mat_u` and `Duu_mat`).

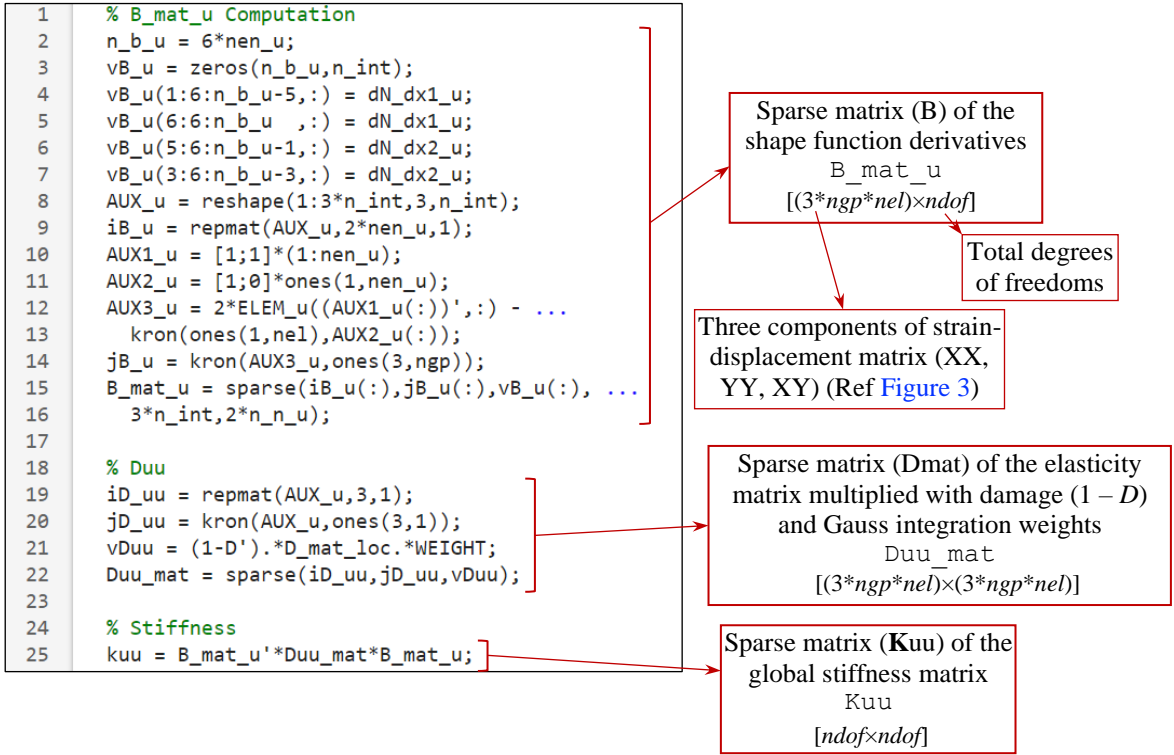

**Figure 7**: Code snippets showing vectorized computation of global stiffness matrix (`Kuu`)

(b) *Variables Update*

The vectorized variable update is shown through code snippets in Figure 8. It is observed in Algorithm 2 (lines 9–16) and Figure 8 that the for-loops over elements and GPs are eliminated. Note that strains, equivalent strains and micro equivalent strains at all the GPs (in the entire model) are evaluated in single lines of codes (Figure 8b, lines 7, 11 and 19) using direct matrix multiplication and entry-wise operations. Moreover, the if-else conditions (used in the non-vectorized code, Figure 5b) are replaced with logical vectors for updating history variables and damage. These logical vectors are called condition vectors (`cond1` and `cond2` in Figure 8b, lines 22 and 27). The fulfillment or non-fulfillment of the conditions are encoded using the not operator (~) in MATLAB. In other words, `cond1` (Figure 8b, line 23) means the condition is fulfilled and `~cond1` (Figure 8b, line 24) means otherwise. The condition vectors are highly efficient in storing the conditions through a single line vector operation that otherwise would need if-else conditions within for-loops.

```
1    % Update primary variables
2    u = u + du;
3
4    % Update SDVs
5    udof = u(1:ndof_u);
6    edof = u(ndof_u + 1:ndof_u + ndof_e);
7    [strain, loc_e_eqv, gp_e_eqv, k_gp, D] = ...
8        update_variables(udof, edof, k0, ki, B_mat_u, ...
9        N_mat_e, alpha_var, beta, c1, c2, c3, k);
```
→ For-Loops over Elements Eliminated

**(a)** Update variables

```
1    % Update Variables......................................
2    function [strain, loc_e_eqv, gp_e_eqv, k_gp, D] = ...
3        update_variables(udof, edof, k0, ki, B_mat_u, ...
4        N_mat_e, alpha_var, beta, c1, c2, c3, k)
5
6    % Local Strain
7    strain = B_mat_u*udof;
8    strain_loc = reshape(strain,3,[])';
9
10   % Local equivalent strain
11   loc_e_eqv = ...
12   c1*(strain_loc(:,1) + strain_loc(:,2))+(1/(2*k))*...
13   ((c2*(strain_loc(:,1) + strain_loc(:,2))).^2) + ...
14   (c3*(strain_loc(:,1).^2 + strain_loc(:,2).^2 - ...
15   strain_loc(:,1).*strain_loc(:,2) + ...
16   3*(strain_loc(:,3).^2)))).^0.5;
17
18   % Equivalent strain
19   gp_e_eqv = N_mat_e*edof;
20
21   % History Equivalent Strain
22   cond1 = (gp_e_eqv - k0) > 1D-10;
23   k_gp(cond1,1) = gp_e_eqv(cond1,1);
24   k_gp(~cond1,1) = k0(~cond1,1);
25
26   % Damage
27   cond2 = (k_gp - ki) < 1D-10;
28   D(cond2,1) = 0.0;
29   D(~cond2,1) = 1-(ki(~cond2,1)./k_gp(~cond2,1)).*...
30      (1 - alpha_var + alpha_var * ...
31      exp(-beta*(k_gp(~cond2,1) - ki(~cond2,1))));
32   end
```
→ For-Loops over GPs Eliminated

→ If-Else conditions eliminated with vectorized conditions

**(b)** Update variables function subroutine

**Figure 8**: Code snippets of the vectorized variable update and its function subroutine

## 2.3 Summary

The MATLAB implementations of the localizing gradient damage method discussed above are summarized through flowcharts in Figure 9. Both non-vectorized and vectorized implementations are shown, along with the parts of code that are vectorized (shown as dashed boxes in Figure 9). The vectorization (shown as dotted lines with arrows in Figure 9) is carried out on the non-vectorized MATLAB implementation to improve its computational efficiency. From Figure 9, it is clear that vectorization reduces the number of operations needed for simulation, consequently increasing the computational speed. Moreover, it can be observed that for-loops needed in the non-vectorized code are entirely eliminated for elemental computations, assembly and variable updates. The for-loops for loadsteps and iterations are unavoidable in a non-linear simulation. Apart from computational efficiency, a vectorized MATLAB code is more compact than a non-vectorized code (Figure 9).

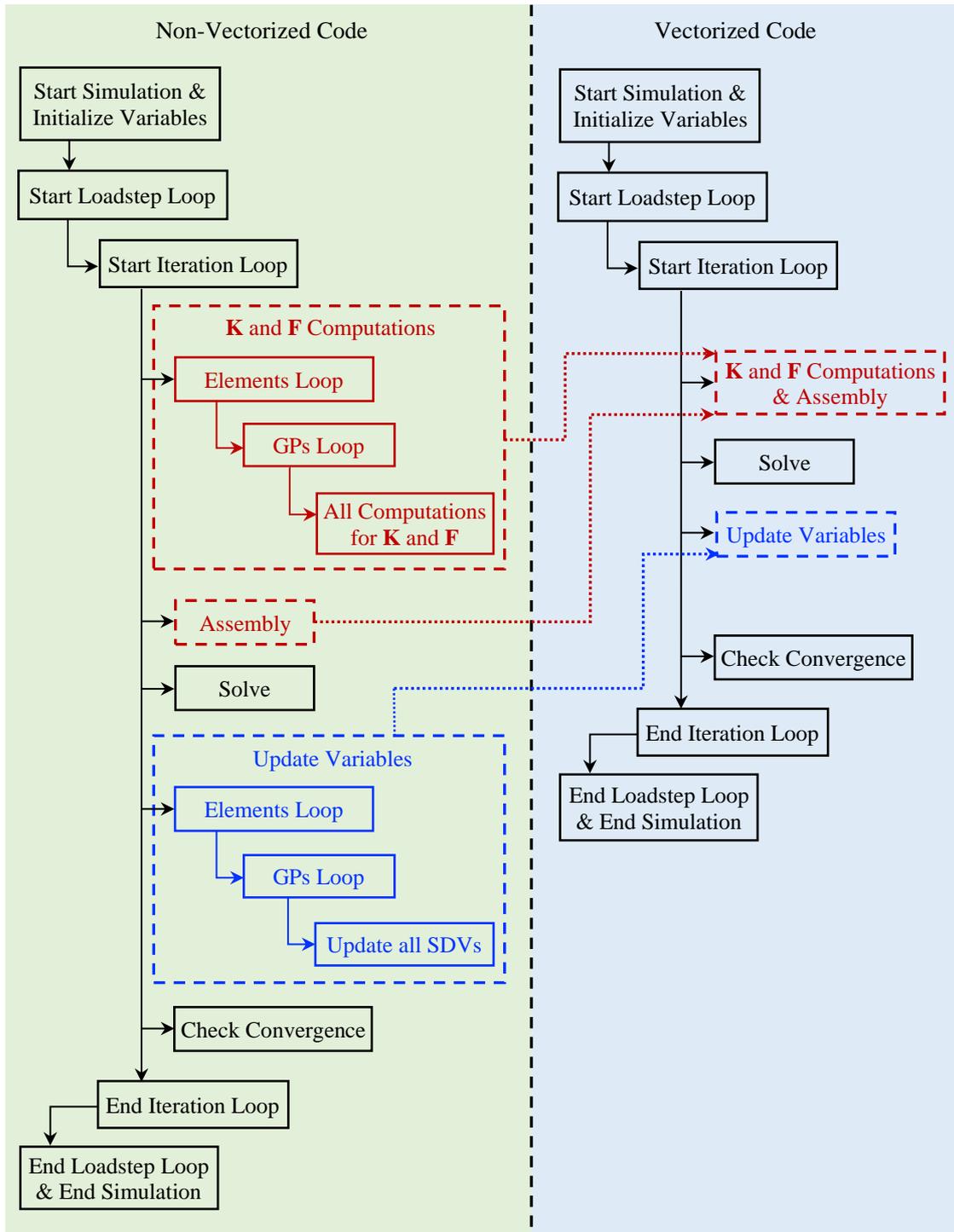

**Figure 9**: Flowchart showing the non-vectorized and vectorized MATLAB implementations for the localizing gradient damage method (LGDM)

## 3. RESULTS AND DISCUSSION

This section compares the non-vectorized and vectorized MATLAB implementations through numerical problems. The comparison is focused on computational efficiency in terms of computation resources (CPU and RAM usage) and computation time. For brevity and straightforward comparison, benchmark numerical problems in 1D, 2D and 3D are chosen whose results are widely available in the published literature. It is emphasized that the proposed MATLAB vectorization is aimed at making LGDM (a non-linear coupled field method) computationally feasible on smaller systems with limited resources. Hence, all the numerical problems in this work are tested with MATLAB 2022a on a consumer-grade PC with specifications – *Processor*: AMD Ryzen 5 5600H; *RAM*: 16 GB; *SSD*: 256 GB.

### 3.1 1D Bar Problem

In this problem, a one-dimensional (1D) bar subjected to tensile loading is simulated using the 1D MATLAB implementation of localizing gradient damage method (LGDM). The problem description and parameters used in the simulation (adopted from Sarkar et al., 2019) are shown in Figure 10. A displacement-controlled load of 0.02 mm is applied on the right side of the bar in 1000 load steps. A defect region is introduced at the middle of the bar to initiate damage. For the simulation, quadratic 1D elements are used for displacement ($u$) and linear 1D elements are used for micro-equivalent strain ($\bar{\varepsilon}_{eq}$).

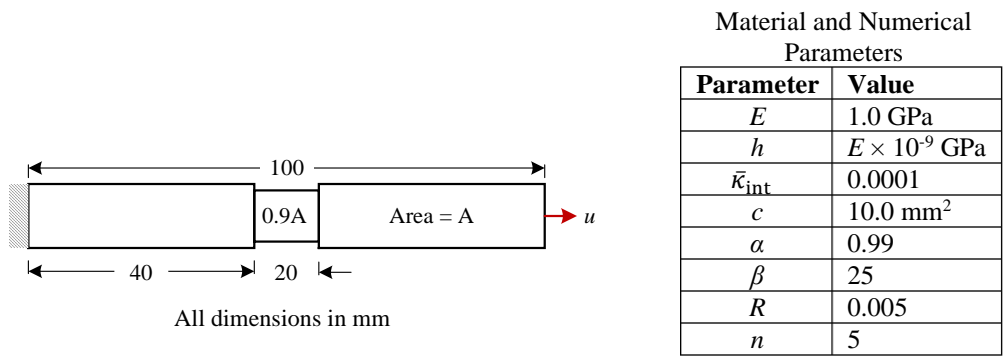

Material and Numerical Parameters

| Parameter | Value |
|---|---|
| $E$ | 1.0 GPa |
| $h$ | $E \times 10^{-9}$ GPa |
| $\bar{\kappa}_{int}$ | 0.0001 |
| $c$ | 10.0 mm$^2$ |
| $\alpha$ | 0.99 |
| $\beta$ | 25 |
| $R$ | 0.005 |
| $n$ | 5 |

**Figure 10**: A schematic representation of the geometry, loads, boundary conditions and parameters for the 1D bar problem

The obtained load-displacement plots for the 1D bar problem are shown in Figure 11. Figure 11a shows a convergence of load-displacement curves obtained using the vectorized implementation using uniform meshes with 500, 800 and 1000 1D elements. It is also shown that the converged load-displacement plot (of 1000 elements) agrees with the reference results available in Sarkar et al. (2019). A comparison of the load-displacement plots (for 1000 elements) obtained using vectorized and non-vectorized implementation in Figure 11b shows that both implementations are equivalent. Apart from these, a similar convergence and equivalence of damage (Figures 12a, b) and micro-equivalent strain (Figures 12c, d) indicate that the vectorized 1D MATLAB implementation is accurate and equivalent to the non-vectorized 1D MATLAB implementation. The obtained plots of damage and micro-equivalent strain are similar to the reference plots in Sarkar *et al.* (2019).

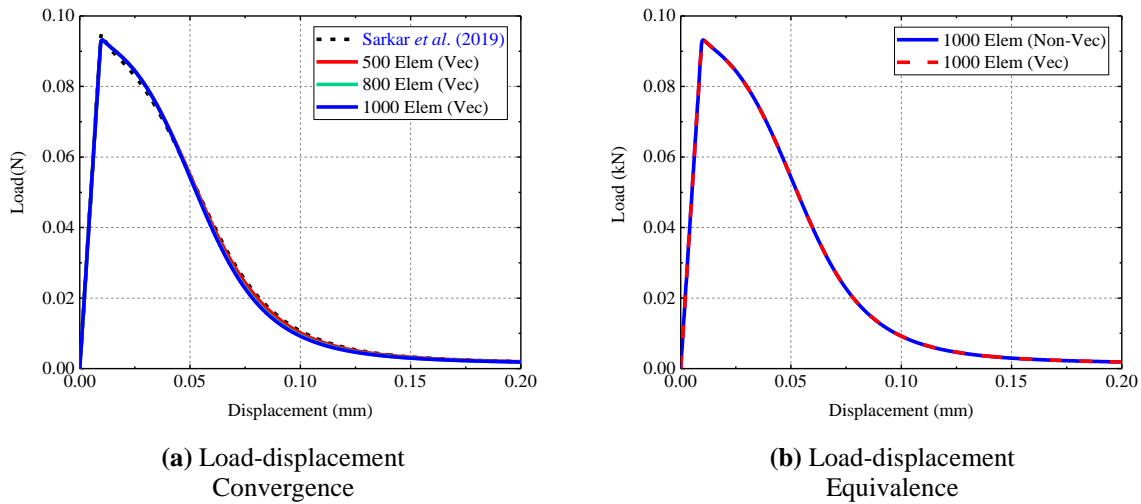

**(a)** Load-displacement Convergence

**(b)** Load-displacement Equivalence

**Figure 11**: Plots showing **(a)** convergence of load-displacement curves obtained using vectorized 1D MATLAB implementation and **(b)** equivalence of load-displacement curves obtained using both vectorized and non-vectorized 1D MATLAB implementation in 1D bar problem

It is established from the equivalence plots that the vectorized and non-vectorized 1D MATLAB implementations yield identical results. However, these implementations are significantly different in terms of computational efficiency. The computational efficiency is compared using computational parameters, i.e. computational resources (CPU and RAM usage) and computation time (wall clock time) taken in the entire simulation. The CPU and RAM usages are calculated by averaging the usages during the entire simulation, and the time taken from the start to the last loadstep is called computation time. The computational parameters are calculated by averaging the values obtained from ten repeated simulations for each case. During these simulations, the MATLAB program is run alone without any other programs/applications simultaneously running on the PC to get the maximum available performance. Typically, a simulation can be regarded as computationally more efficient when it fulfills the following three criteria,

   (a) CPU usage – High (↑): Implying usage of multiple processor cores
   (b) RAM usage – Low (↓): Implying efficient/less memory handling
   (c) Computation time – Low (↓): Implying the effect of the above two factors

Another parameter, called time per iteration, is used to show the contribution of different code sections to the overall computation time for each iteration. This parameter identifies the parts of code that are either more efficient or inefficient. The time per iteration is calculated for three major parts of the code, which are,

(a) Assembly:        Algorithm 1, lines 4 – 11;    Algorithm 2, lines 6 – 7
(b) Solve:           Algorithm 1, line 12;          Algorithm 2, line 8
(c) Variable update: Algorithm 1, lines 13 – 24;    Algorithm 2, lines 9 – 16

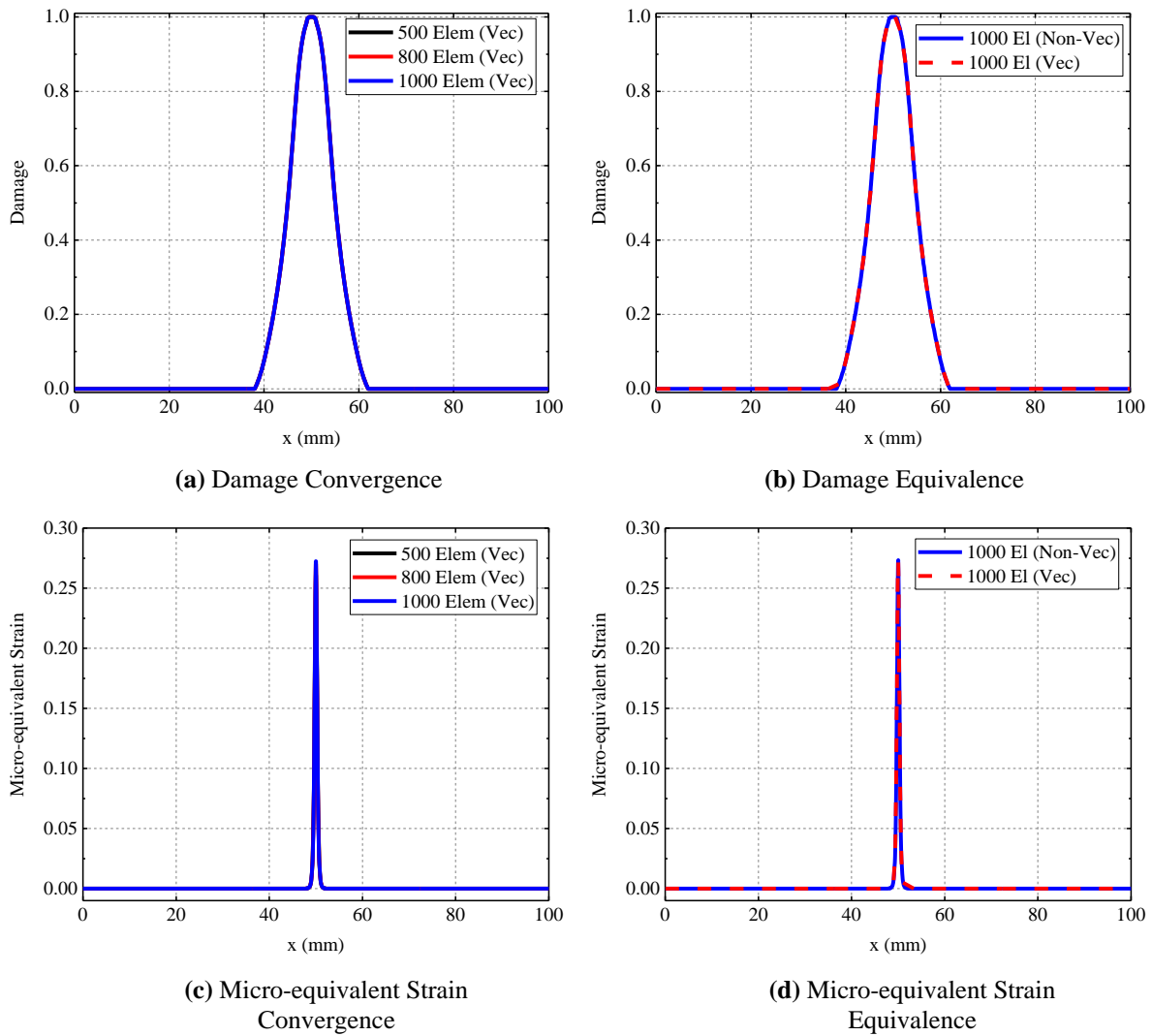

**Figure 12**: Plots showing convergence of **(a)** damage; **(c)** micro-equivalent strain obtained using vectorized 1D MATLAB implementation and equivalence of **(b)** damage; **(d)** micro-equivalent strain obtained using both vectorized and non-vectorized 1D MATLAB implementation in 1D bar problem

The comparison of computational parameters and time per iteration obtained from vectorized and non-vectorized 1D MATLAB implementation is shown in Figure 13. The computational parameters in Figure 13a are obtained by simulating the 1D problem using 1000 elements. In Figure 13b, the individual operation times (assembly, solve and variable update) are obtained by dividing the total time (for that operation in the entire simulation) by the total number of iterations in the entire simulation.

The computational efficiency of the vectorized implementation is evident in Figure 13. It is observed that despite a marginal change in CPU (~ 1% ↑) and RAM usage (~0.08 GB ↓), the computation time (~90% ↓) is significantly lower for the vectorized implementation. The low computation time is also reflected in the single iteration times shown in Figure 13b. It can be observed that significant reductions in the assembly and variable update times resulted in bringing down the total computation time. The change in solve time is marginal due to the same direct solver (`mldivide`) being used for both implementations. The marginal improvements in CPU and RAM usage may be attributed to almost equivalent memory handling in the simulation and serial execution. However, the performance differences are more pronounced in the 2D and 3D problems discussed later in the following sub-sections.

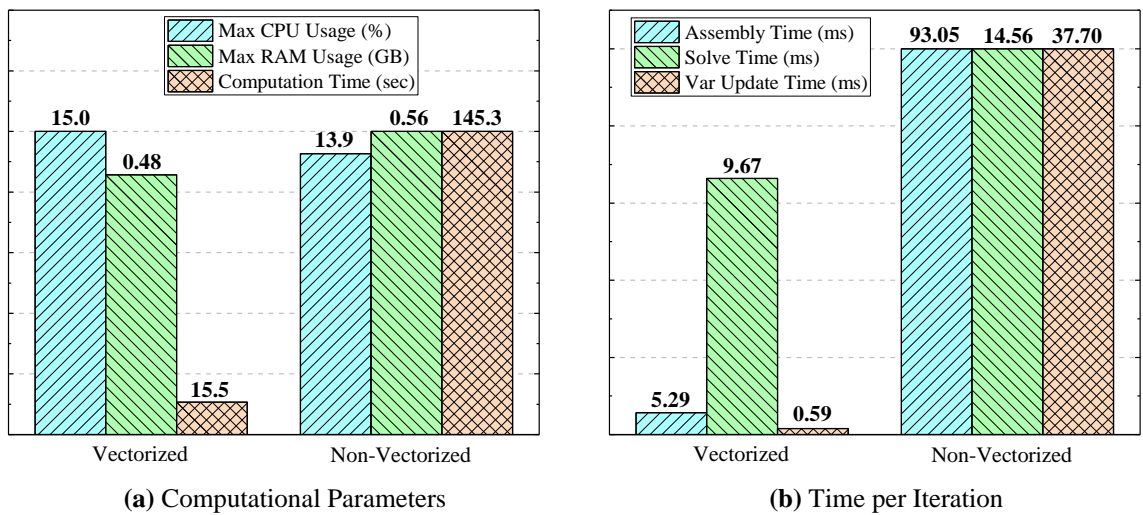

**(a)** Computational Parameters  **(b)** Time per Iteration

**Figure 13**: Plots showing comparison of **(a)** Computational parameters and **(b)** Time per iteration obtained using both vectorized and non-vectorized 1D MATLAB implementation in 1D bar problem

Based on the abovementioned observations and criteria defined previously, the vectorized MATLAB implementation in 1D can be considered computationally more efficient than the non-vectorized implementation. The outcome of the comparison is summarized in Figure 18. In Figure 18, the change of a computational parameter is considered favorable when it fulfills the criteria of computational efficiency. The changes in the CPU and RAM usage, although marginal, are favorable i.e. CPU usage increased and RAM usage decreased. However, the change in computation time is significant and also favorable.

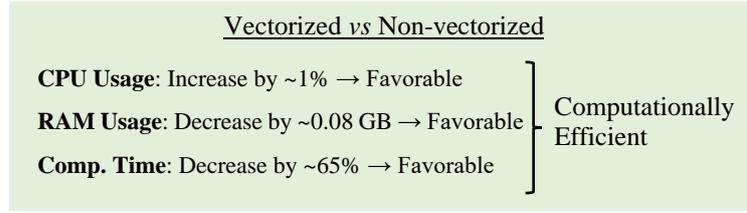

Vectorized *vs* Non-vectorized

**CPU Usage**: Increase by ~1% → Favorable
**RAM Usage**: Decrease by ~0.08 GB → Favorable    } Computationally Efficient
**Comp. Time**: Decrease by ~65% → Favorable

**Figure 14**: Outcome of comparisons between the vectorized and non-vectorized 1D MATLAB implementations

### 3.2 2D Side Edge Notch (SEN) Problem

In this problem, a side edge notch specimen is simulated using the 2D MATLAB implementation of the localizing gradient damage method (LGDM). The problem description and parameters used in the simulation (adopted from Sarkar *et al.*, 2019) are shown in Figure 15. A displacement-controlled load of 0.8 mm is applied on the top surface in 80 load steps of 0.01 mm. For the simulation, uniform finite element meshes with 50×50 (2500), 100×100 (10000) and 120×120 (14400) elements are used. The model uses 2D biquadratic elements for displacement ($u$) and 2D bilinear elements for micro-equivalent strain ($\bar{\varepsilon}_{eq}$).

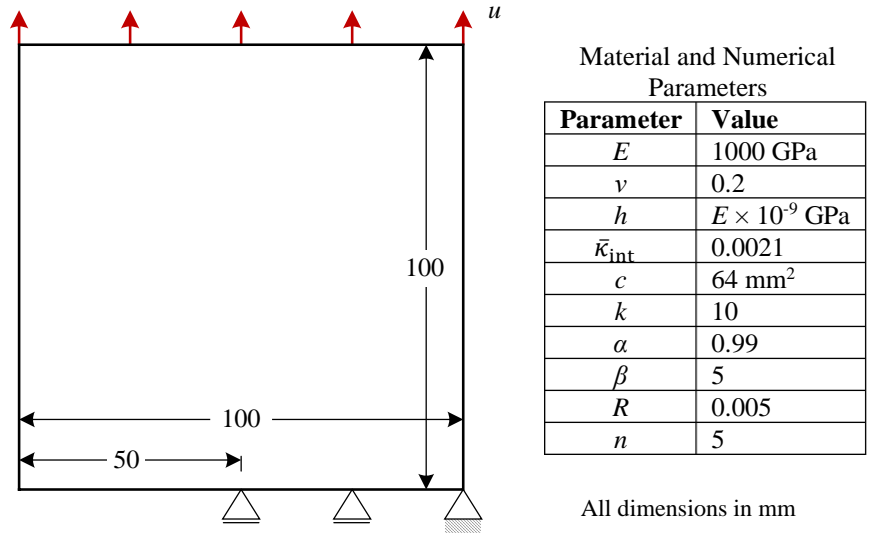

**Figure 15**: A schematic representation of the geometry, loads, boundary conditions and parameters for the 2D side edge notch (SEN) problem

Similar to the 1D bar problem, a convergence of load-displacement curves and equivalence of non-vectorized and vectorized 2D MATLAB implementations is shown in Figure 16. For validation, it is observed in Figure 16a that the converged load-displacement curve (for the 100×100 elements mesh) obtained from the vectorized implementation is almost similar to the reference results in Sarkar *et al.* (2019). Figure 16b shows the equivalence of load-displacement curves obtained from non-vectorized and vectorized implementations. Additionally, the damage and micro-equivalent strain plots obtained from vectorized implementation using 100×100 elements mesh are shown in Figure 17. The obtained damage and micro-equivalent strain plots are found similar to the reference plots in Sarkar *et al.* (2019). Hence, they establish the accuracy and consistency of the vectorized MATLAB implementation in 2D.

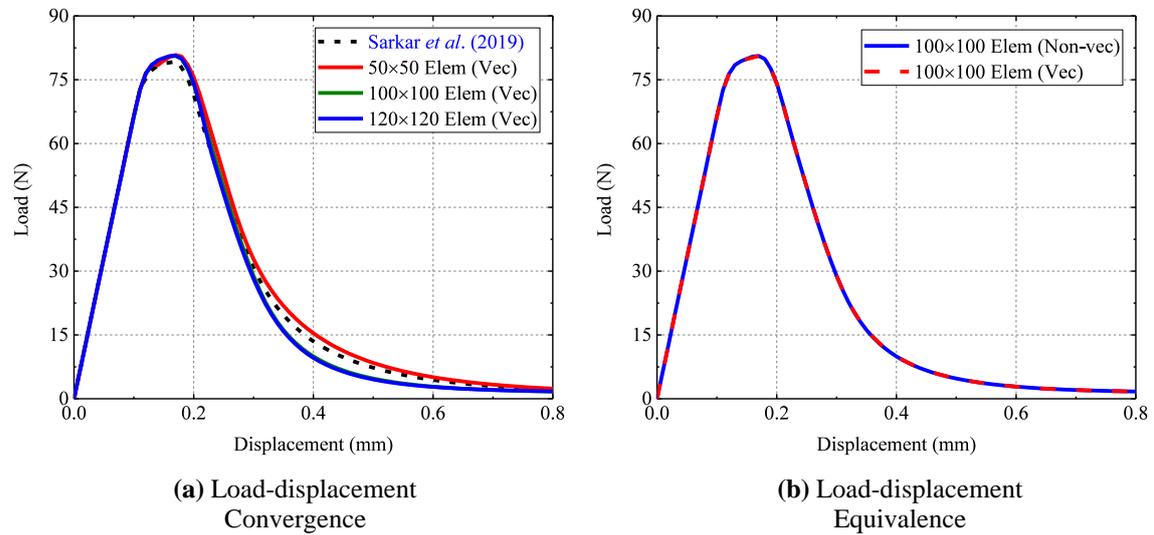

**(a)** Load-displacement Convergence

**(b)** Load-displacement Equivalence

**Figure 16**: Plots showing **(a)** convergence of load-displacement curves obtained using vectorized 2D MATLAB implementation and **(b)** equivalence of load-displacement curves obtained using both vectorized and non-vectorized 2D MATLAB implementation in 2D SEN problem

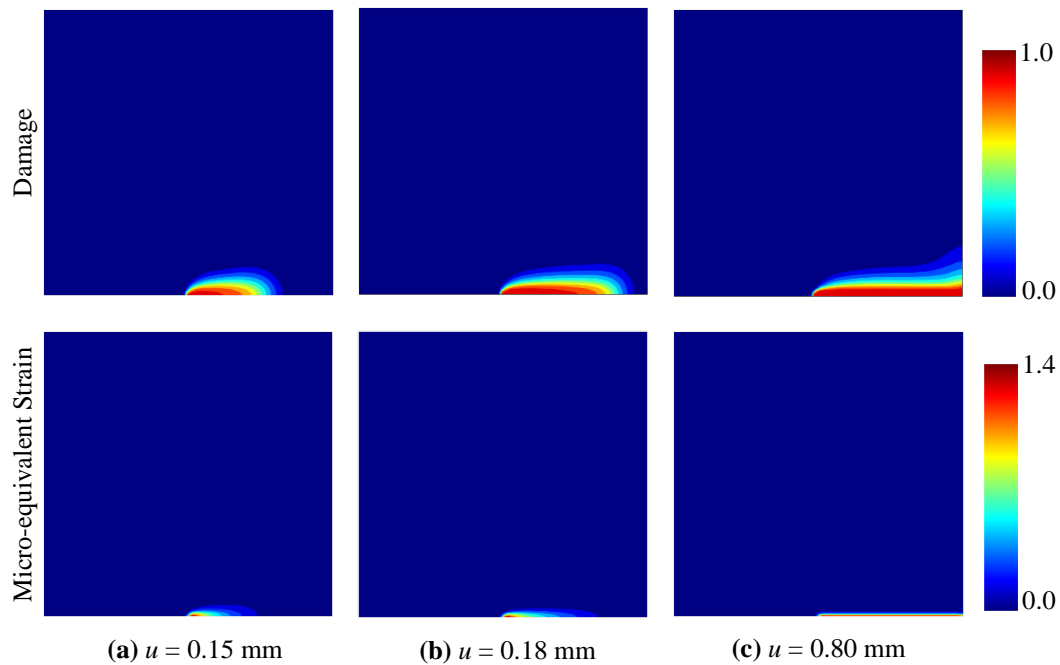

**(a)** $u$ = 0.15 mm  **(b)** $u$ = 0.18 mm  **(c)** $u$ = 0.80 mm

**Figure 17**: Plots showing evolution of damage and micro-equivalent strain obtained at different applied displacements using vectorized 2D MATLAB implementation in 2D SEN problem

To this end, it can be stated that both non-vectorized and vectorized 2D MATLAB implementations are accurate and equivalent. However, their computational efficiencies are found to be significantly different. Similar to the 1D problem in the previous sub-section, computational parameters and time per iteration obtained from the 2D SEN problem using the 100×100 elements mesh are shown in Figure 18. It is noted from Figure 18 that two types of non-vectorized implementation are tested, namely, with and without `parfor` (`parfor` is a MATLAB feature that parallelizes for-loops). In this work, `parfor` is used for 2D and 3D simulations only because it led to decreased computational efficiency in 1D simulations.

From Figure 18, the comparison of vectorized implementation with the non-vectorized implementation (with and without `parfor`) is carried out as follows:

- *Vectorized vs. non-vectorized implementation (**without** `parfor`):*

  The computational parameters of non-vectorized implementation (without `parfor`) are shown in Figure 18a (middle plot). It is observed that the CPU usage of vectorized implementation (Figure 18a, left plot) is almost double (~100% ↑) of the non-vectorized implementation, possibly due to parallel/vectorized computations. In contrast, the RAM usage of vectorized implementation is slightly less (~0.2 GB ↓) compared to the non-vectorized implementation. This can be attributed to the serial execution of both vectorized and non-vectorized (without `parfor`) implementations resulting in almost unchanged memory usage.

  Nonetheless, higher CPU usage during the simulation might have led to a significantly less computation time (~65% ↓) of the vectorized implementation compared to the non-vectorized implementation (without `parfor`). In Figure 18b, a similar trend with significantly less assembly and variable update times (per iteration) for the vectorized implementation is observed compared to non-vectorized implementation (without `parfor`). The reduction in computation time during assembly and variable updates contributes significantly to the overall decrease in the computation time for the vectorized implementation.

- *Vectorized vs. non-vectorized implementation (**with** `parfor`):*

  It is observed that the use of MATLAB `parfor` in non-vectorized implementation (Figure 18a, right plot) leads to a significant increase in CPU usage and decrease computation time compared to the non-vectorized implementation without `parfor` (Figure 18a, middle plot). However, using `parfor` increases the RAM usage of the non-vectorized implementation (with `parfor`) by more than double compared to the non-vectorized implementation (without `parfor`). The RAM usage is even more when compared with the vectorized implementation that has ~2.8 GB less usage (Figure 17a, left plot). This makes the non-vectorized implementation (with `parfor`) computationally inefficient compared to the vectorized implementation. The high RAM usage can be attributed to increased memory needed for the execution of `parfor` that distributes data and processes of the for-loops across multiple CPU processors (MATLAB, 2022).

Apart from less RAM usage, it is observed (in Figure 18a) that the computation time of the vectorized implementation is considerably less (~40% ↓) compared to the non-vectorized implementation (with `parfor`) despite unfavorable CPU usage (~32% ↓). In other words, decreased CPU usage usually implies fewer computations and slower simulation, thus, unfavorable. However, decreased CPU usage in addition to less computation time may imply that the vectorized implementation needs fewer computations than the non-vectorized implementation. This makes the vectorized implementation computationally more efficient than the non-vectorized implementation (with `parfor`). Similar trends of less assembly and variable update times are observed for vectorized implementation in the time per iteration plots shown in Figure 18b.

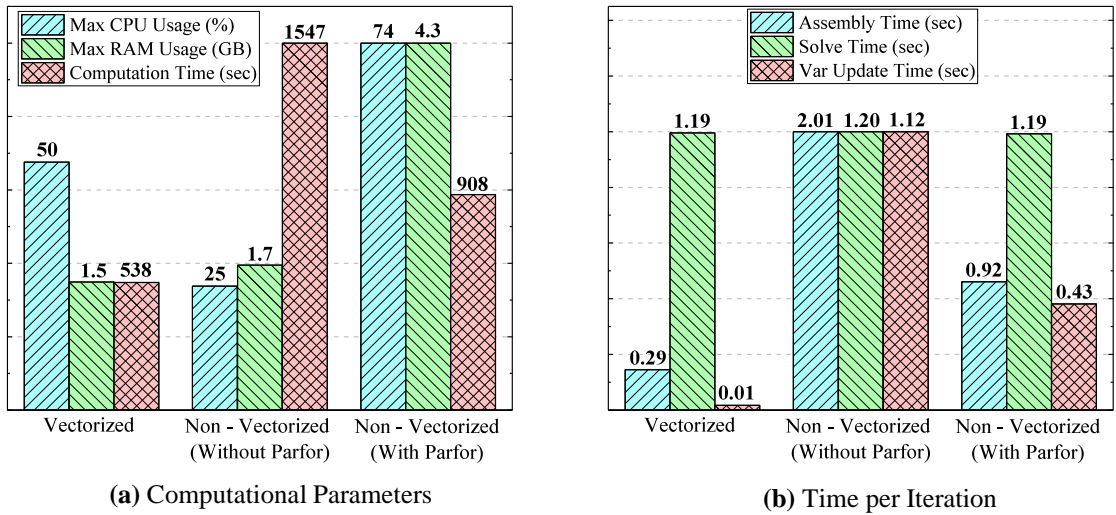

**(a)** Computational Parameters    **(b)** Time per Iteration

**Figure 18**: Plots showing comparison of **(a)** Computational parameters and **(b)** Time per iteration obtained using both vectorized and non-vectorized 2D MATLAB implementation in 2D SEN problem

- *Comparison Summary*:

  The abovementioned observations and comparisons indicate that using MATLAB `parfor` makes the non-vectorized implementation computationally more efficient in 2D simulations. However, even with `parfor`, non-vectorized implementation cannot beat the computation time and RAM usage of vectorized implementation. Moreover, the higher RAM usage by MATLAB `parfor` limits its use for simulating larger models on limited resources. The computational efficiency of the 2D vectorized MATLAB implementation over the non-vectorized implementation (with and without `parfor`) is summarized in Figure 19.

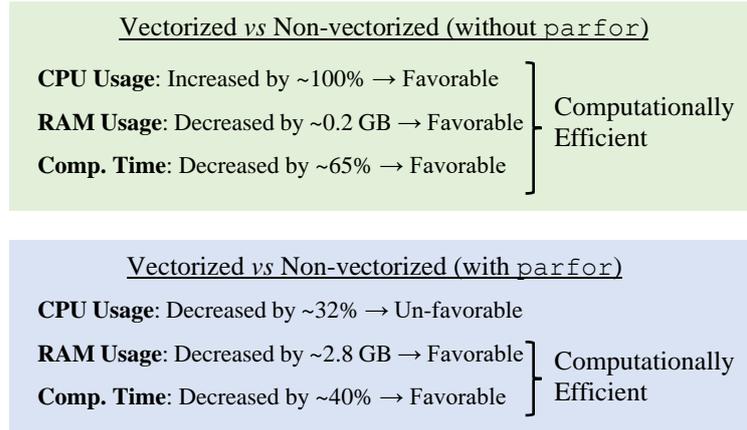

**Figure 19**: Outcome of comparisons between the vectorized and non-vectorized 2D MATLAB implementations

## 3.3 3D Side Edge Notch (SEN) Problem

In this problem, a 3D side edge notch (SEN) specimen is simulated using the 3D MATLAB implementation of the localizing gradient damage method (LGDM). This problem is an extension of the 2D SEN problem (solved in the previous sub-section) with a thickness of 10 mm in the third dimension. The problem description and parameters used in the simulation (adopted from Sarkar *et al*., 2019) are shown in Figure 20. A displacement-controlled load of 0.8 mm is applied on the top surface in 80 load steps of 0.01 mm. For the simulation, uniform finite element meshes with 50×50×5 (12500), 80×80×5 (32000) and 100×100×5 (50000) elements are used. 3D trilinear elements are used for both displacement (*u*) and micro-equivalent strain ($\bar{\varepsilon}_{eq}$). Note that trilinear elements (instead of triqudratic) are used for displacement (*u*) primarily to reduce computational effort on the consumer-grade PC used for simulations. The use of trilinear elements (for *u*) is not a limitation of the current 3D MATLAB implementation, which can be easily modified to include triquadratic elements. Moreover, it is pointed out that using trilinear elements instead of triquadratic elements in the LGDM leads to similar results without any compromise in accuracy, as shown in Sarkar *et al*. (2021a) and Sarkar *et al*. (2022b).

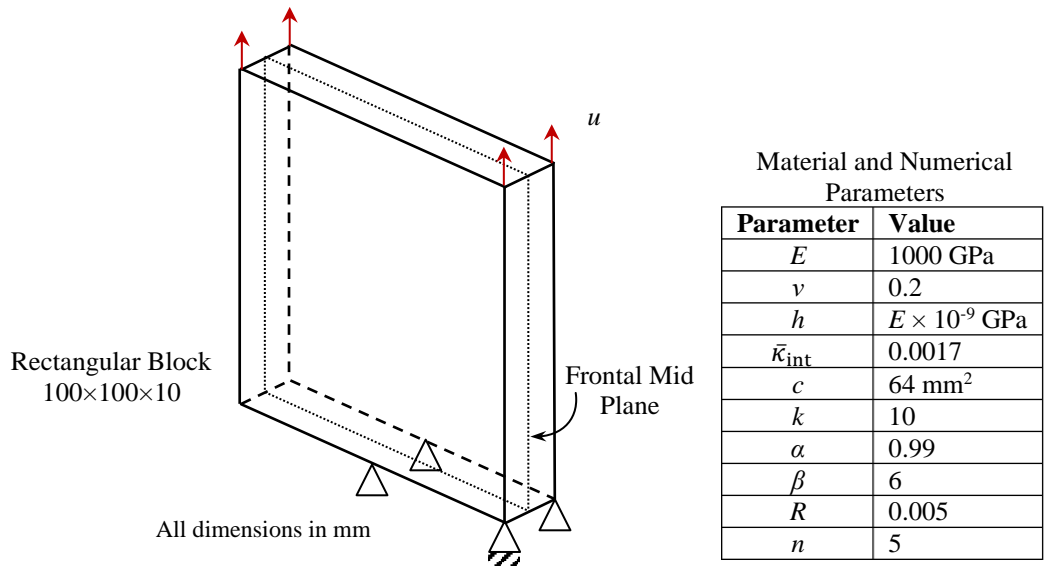

**Figure 20**: A schematic representation of the geometry, loads, boundary conditions and parameters for the 3D SEN problem

Similar to the previous problems, the convergence of load-displacement curves and equivalence of non-vectorized and vectorized 3D MATLAB implementations is shown in Figure 21. For validation, it is shown in Figure 21a that the load-displacement curves obtained from the vectorized implementation converge (at the mesh of 100×100×5 elements) and are similar to the reference results in Sarkar *et al*. (2019). The reference results reported by Sarkar *et al*. (2019) are for a 2D plane strain case (with unit thickness). Therefore, the load-displacement curves from the present problem (3D simulation) are scaled down according to the thickness for comparison with the reference results. Figure 20b shows the equivalence plot of load-displacement curves obtained from both vectorized and non-vectorized implementations. The equivalence plot shows that the results of vectorized and non-vectorized 3D MATLAB implementations are identical.

Apart from these, the damage and micro-equivalent strain plots obtained from vectorized implementation using 100×100×5 elements mesh are shown in Figures 22-23. The obtained damage and micro-equivalent strain plots are similar to the reference 2D simulation plots in Sarkar *et al*. (2019), hence establishing the accuracy of the vectorized 3D MATLAB implementation. The damage evolution is also shown as iso-surface plots in Figure 24, where the typical thumbnail shape of damage evolution is observed (Yamamoto *et al*., 1987). The 3D regions (shown in red) in iso-surface plots of Figure 24 indicate damage ($D$) > 0.9.

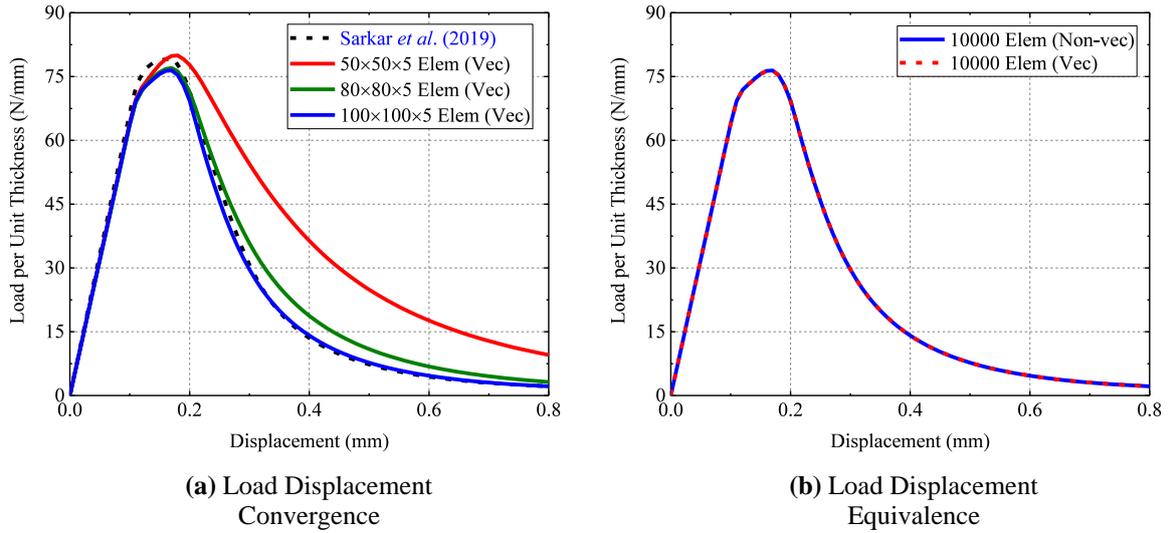

**(a)** Load Displacement Convergence

**(b)** Load Displacement Equivalence

**Figure 21**: Plots showing **(a)** Convergence of load-displacement curves obtained using vectorized 3D MATLAB implementation and **(b)** Equivalence of load-displacement curves obtained using both vectorized and non-vectorized 3D MATLAB implementation in 3D SEN problem

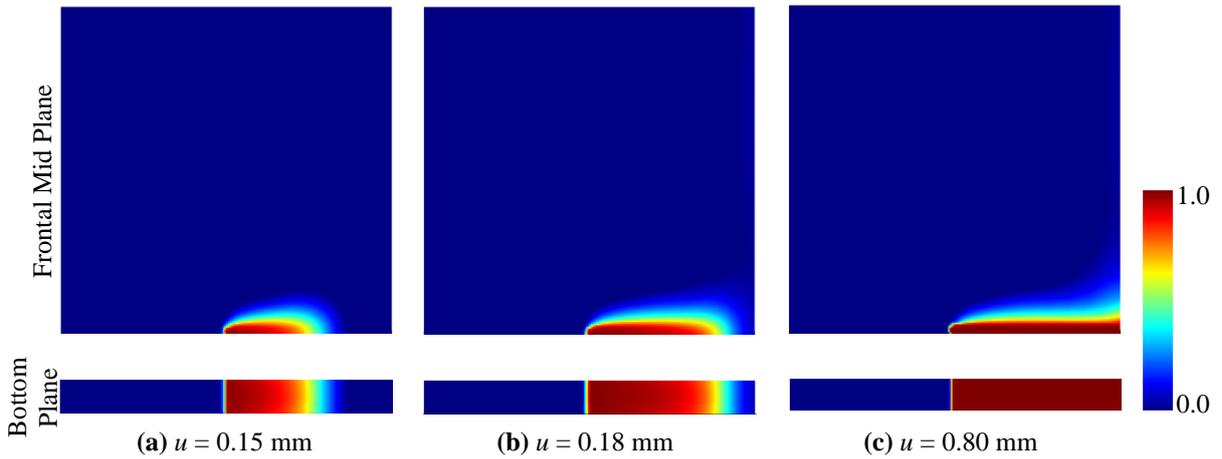

**(a)** $u = 0.15$ mm   **(b)** $u = 0.18$ mm   **(c)** $u = 0.80$ mm

**Figure 22**: Plots showing evolution of damage obtained at different applied displacements using vectorized 3D MATLAB implementation in 3D SEN problem

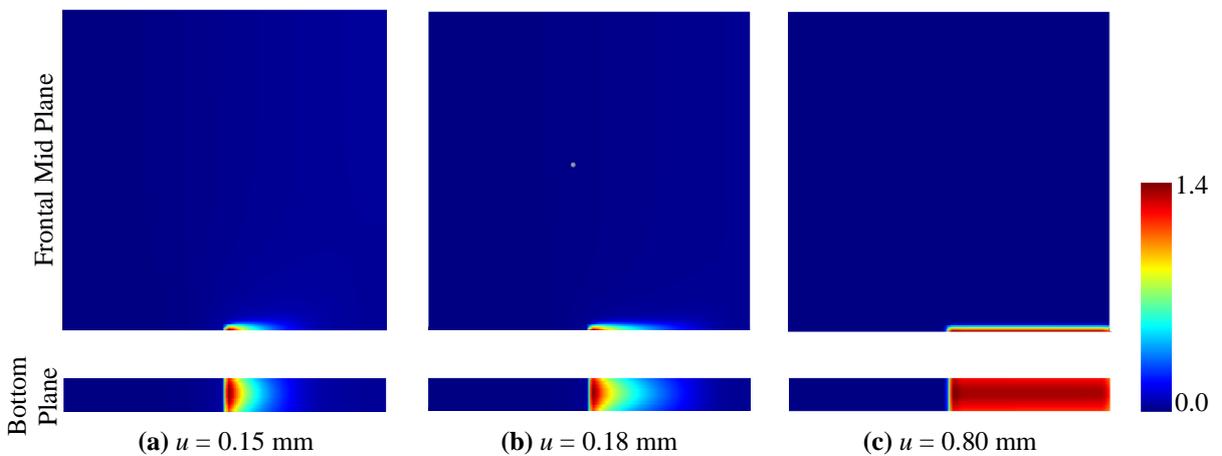

**(a)** $u = 0.15$ mm   **(b)** $u = 0.18$ mm   **(c)** $u = 0.80$ mm

**Figure 23**: Plots showing evolution of micro-equivalent strain obtained at different applied displacements using vectorized 3D MATLAB implementation in 3D SEN problem

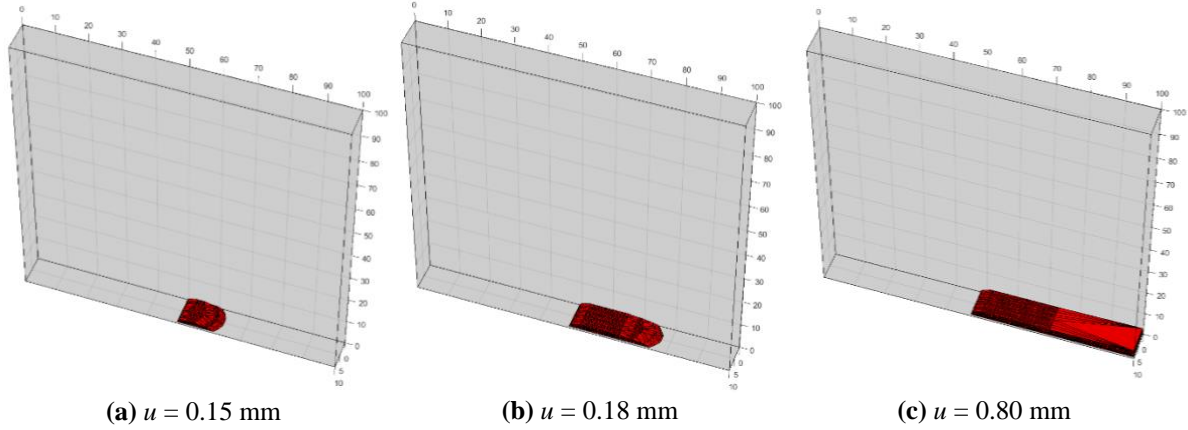

**(a)** *u* = 0.15 mm      **(b)** *u* = 0.18 mm      **(c)** *u* = 0.80 mm

**Figure 24**: Iso-surface plots showing evolution of damage region ($D > 0.9$) obtained at different applied displacements using vectorized 3D MATLAB implementation in 3D SEN problem

Similar to the 1D and 2D problems in the previous sub-sections, computational parameters and time per iteration obtained from the 3D SEN problem using the 100×100×5 elements mesh are shown in Figure 25. The vectorized implementation is compared with the non-vectorized implementation (with and without `parfor`) as follows,

- *Vectorized vs. non-vectorized implementation* (***without*** `parfor`):

  The computational parameters of non-vectorized implementation (without `parfor`) are shown in the middle plot of Figure 25a. It is observed that the computation time of vectorized implementation (Figure 25a, left plot) is less than half (~56% ↓) of the non-vectorized implementation (without `parfor`) while the RAM usage is almost similar (~0.1 GB ↓). The similar RAM usage can be attributed to the serial execution of both vectorized and non-vectorized (without `parfor`) implementations resulting in almost unchanged memory usage. Additionally, the CPU usage of the vectorized implementation is less by ~30% (↓), possibly due to fewer computations, resulting in less computation time.

  In Figure 25b, a similar trend with a significant decrease in assembly and variable update times (per iteration) for the vectorized implementation is observed compared to non-vectorized implementation (without `parfor`). The reduction in computation time during assembly and variable updates contributes significantly to the overall decrease in the computation time.

- *Vectorized vs. non-vectorized implementation* (***with*** `parfor`):

  It is found that the non-vectorized implementation (Figure 25a, right plot) using `parfor` is more efficient than the non-vectorized implementation without `parfor` (Figure 25a, middle plot). The efficiency is evident from the fact that the non-vectorized implementation (with `parfor`) has a computation time of almost half of the non-vectorized implementation (without `parfor`). However, the RAM usage in non-vectorized implementation (with `parfor`) increases to almost twice due to MATLAB `parfor` while increasing the CPU usage slightly.

Compared to the non-vectorized implementation (with `parfor`), the vectorized implementation (Figure 25a, left plot) has less computation time (~11% ↓). Additionally, the RAM (~3.8 GB ↓) and CPU (~40% ↓) usage of the vectorized implementation are significantly less. The decreased CPU usage may imply slower computations, but a decrease in computation time (along with CPU usage) suggests less number of computations during the simulations. Hence, the vectorized implementation can be considered computationally more efficient than the non-vectorized implementation (with `parfor`).

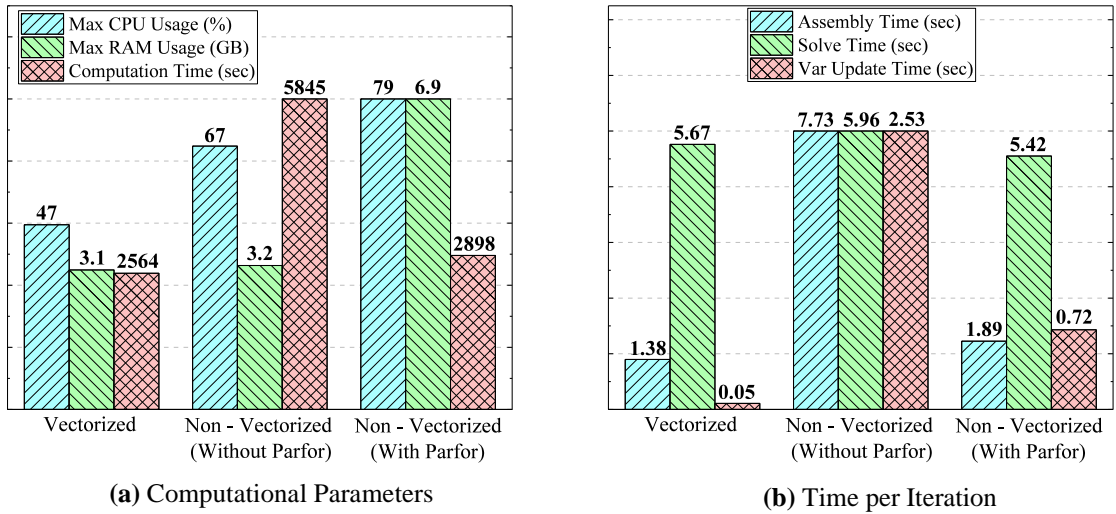

(a) Computational Parameters  (b) Time per Iteration

**Figure 25**: Plots showing comparison of **(a)** Computational parameters and **(b)** Time per iteration obtained using both vectorized and non-vectorized 3D MATLAB implementation in 3D SEN problem

- *Comparison Summary*:

  The abovementioned observations and comparisons indicate that using MATLAB `parfor` makes the non-vectorized implementation computationally more efficient in 3D simulations. However, the computational cost (in terms of RAM usage) of running MATLAB `parfor` in 3D is significantly higher than 1D or 2D simulations. The number of elements in the z-direction is kept at five in the current problem to keep the computation effort within the limits of the consumer-grade PC used for simulations. With a higher number of elements, the computations exceeded the PC's RAM capacity and started using SSD, thus, making simulation slower. Moreover, the RAM capacity is exceeded more frequently in the simulations that used MATLAB `parfor`. It is observed that the vectorized 3D MATLAB implementation takes less computation time along with less RAM and CPU usage, making the simulation of even larger models possible. The simulation of these larger models would have been impossible using the non-vectorized 3D implementation with available computational resources.

Nonetheless, all MATLAB implementations in the present study are designed to simulate models (with a higher number of elements) as efficiently as possible using given computational resources. The computational efficiency of the vectorized 3D MATLAB implementation compared to the non-vectorized implementation (with and without `parfor`) is summarized in Figure 26.

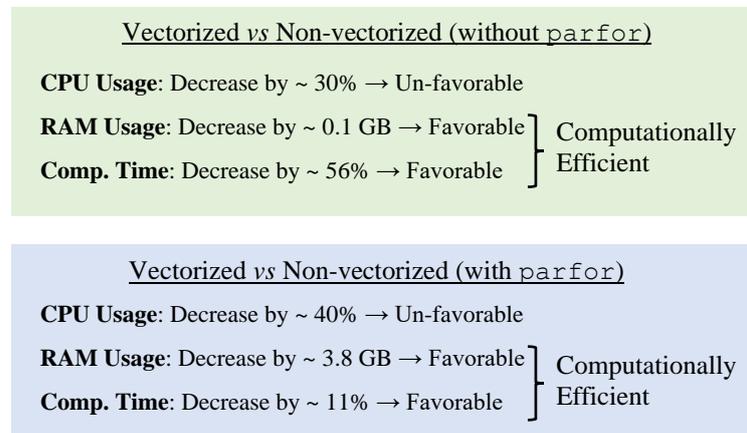

**Figure 26**: Outcome of comparisons between the vectorized and non-vectorized 3D MATLAB implementations

## 4. CONCLUSION

In this work, MATLAB implementations for the localizing gradient damage method (LGDM) are presented, that include both non-vectorized and vectorized MATLAB implementations. It is shown through numerical problems in 1D, 2D and 3D that the vectorized MATLAB implementations are computationally efficient in terms of RAM usage and computation time compared to the non-vectorized MATLAB implementations. Various details of both non-vectorized and vectorized implementations in MATLAB are discussed and their differences are highlighted. The presented discussion along with the provided source codes can be used to easily convert any non-linear finite element based non-vectorized code to vectorized code for increasing computational efficiency.

The major conclusions drawn from the present work are as follows,

- Both vectorized and non-vectorized MATLAB implementations of LGDM are accurate and equivalent in terms of damage simulation and its results
- The use of MATLAB `parfor` in the non-vectorized implementation leads to an increase in computational efficiency of the non-vectorized implementation
- The computational efficiency of the MATLAB implementations in decreasing order is, Vectorized > Non-vectorized (with `parfor`) > Non-vectorized (without `parfor`)
- The higher computational efficiency of vectorized implementation makes it suitable for carrying out simulations on systems with less computational resources and enables the simulation of even larger models with available resources

## APPENDIX A

The damage ($D$) law used in the simulations is shown in Eq. A.1.

$$D(\bar{\kappa}) = \begin{cases} 1 - \frac{\bar{\kappa}_0}{\bar{\kappa}}\{1 - \alpha + \alpha \exp[-\beta(\bar{\kappa} - \bar{\kappa}_0)]\}, & \bar{\kappa} > \bar{\kappa}_0 \\ 0, & \bar{\kappa} \leq \bar{\kappa}_0 \end{cases} \quad (A.1)$$

where, $\bar{\kappa}$ is the history micro-equivalent strain (defined in Eq. A.2), $\alpha$ and $\beta$ are material parameters. In Eq. A.2, $\bar{\varepsilon}_{eq}$ is the value of micro-equivalent strain at an instantaneous time $\tau$ during the entire loading time (0 to $t$).

$$\bar{\kappa}(t) = \max\{\bar{\varepsilon}_{eq} | 0 \leq \tau \leq t\} \quad (A.2)$$

The definition of the equivalent strain called the modified von Mises strain is shown in Eq. A.3. In this equation, $I_1$ and $J_2$ are the invariants of the strain tensor while $\nu$ is the Poisson's ratio and $k$ is the ratio of compressive to tensile strength.

$$\varepsilon_{eq} = \frac{k-1}{2k(1-2\nu)}I_1 + \frac{1}{2k}\sqrt{\frac{(k-1)^2}{(1-2\nu)^2}I_1^2 + \frac{2k}{(1-\nu)^2}J_2} \quad (A.3)$$

The interaction function ($g$) used in the formulation is shown in Eq. A.4, in which, $R$ and $n$ are the material parameters.

$$g(D) = \frac{(1-R)e^{(-nD)} + R - e^{(-n)}}{1 - e^{(-n)}} \quad (A.4)$$

## ACKNOWLEDGEMENT

This research did not receive any specific grant from funding agencies in the public, commercial, or not-for-profit sectors. However, the author is grateful to the Rensselaer

Polytechnic Institute, Troy, NY for providing access to MATLAB software through its academic license.